\documentclass{article}%
\usepackage{amsmath,amsfonts,amssymb,amsthm, bbm}
\usepackage{amsmath}
\usepackage{amsfonts}
\usepackage{amssymb}%
\providecommand{\U}[1]{\protect\rule{.1in}{.1in}}
\numberwithin{equation}{section}
\newtheorem{lemma}{Lemma}[section]
\newtheorem{theorem}{Theorem}[section]
\newtheorem{corollary}{Corollary}[section]
\newtheorem{proposition}{Proposition}[section]
\newtheorem{definition}{Definition}[section]

\newtheorem{counterexample}{Counterxample}[section]
\newtheorem{remark}{Remark}[section]
\begin{document}

\title{\vspace{-1in}\parbox{\linewidth}{\footnotesize\noindent
} \vspace{\bigskipamount} \\Backward stochastic differential equations under \\super linear G-expectation and associated \\Hamilton-Jacobi-Bellman equations \thanks{The author would like to thank Prof.
S. Peng for instructions and Dr. Mingshang Hu for his help.} }
\author{Yuhong Xu \thanks{Institute of Mathematics, Shandong University, Jinan 250100,
P.R.China, e-mail: xuyuhongmath@163.com.}}
\date{}
\maketitle

\indent\textbf{Abstract.} This paper first studies super linear G-expectation.
Uniqueness and existence theorem for backward stochastic differential
equations (BSDEs) under super linear expectation is established to provide
probabilistic interpretation for the viscosity solution of a class of
Hamilton-Jacobi-Bellman equations, including the well known
Black-Scholes-Barrenblett equation, arising in the uncertainty volatility
model in mathematical finance. We also show that BSDEs under super linear
expectation could characterize a class of stochastic control problems. A
direct connection between recursive super (sub) strategies with mutually
singular probability measures and classical stochastic control problems is
provided. By this result we give representation for solutions of
Black-Scholes-Barrenblett equations and G-heat equations. \newline

\indent\textbf{Key words.} G-Brownian motion, super linear expectation, normal
distribution, backward stochastic differential equations,
Hamilton-Jacobi-Bellman equation, Feynman-Kac formula, stochastic
optimization, uncertainty volatility model \newline

\indent\textbf{AMS subject classifications.} 60H10, 60H30, 60J65, 35K55,
35K05, 49L25

\section{Introduction}

The motive of this paper is to show that backward stochastic differential
equations (BSDEs) under super linear expectation coupled with a forward
diffusion:%
\begin{align}
X_{s}^{t,x}  &  =x+\int_{t}^{s}b\left(  X_{r}^{t,x}\right)  dr+\int_{t}%
^{s}h_{j}\left(  X_{r}^{t,x}\right)  d\left\langle B^{j}\right\rangle
_{r}+\int_{t}^{s}\sigma_{j}\left(  X_{r}^{t,x}\right)  dB_{r}^{j}%
,\ t\in\lbrack0,T],\\
Y_{s}^{t,x}  &  =\mathbb{E}_{\ast}[\Phi(X_{T})+\int_{s}^{T}g\left(
X_{r}^{t,x},Y_{r}^{t,x}\right)  dr+\int_{s}^{T}f_{j}\left(  X_{r}^{t,x}%
,Y_{r}^{t,x}\right)  d\left\langle B^{j}\right\rangle _{r}|\mathcal{F}%
_{s}],s\in\lbrack t,T],
\end{align}
provide a probabilistic interpretation for the viscosity solution of a class
of Hamilton-Jacobi-Bellman equations (HJB):%

\begin{align}
\partial_{t}u+\underset{\alpha\in\Gamma}{\inf}\left\{  \mathcal{L}\left(
x,\alpha\right)  u+g\left(  x,u\right)  \right\}   &  =0\\
u|_{t=T}  &  =\Phi.
\end{align}
where $\mathcal{L}\left(  x,\alpha\right)  $ is a second order elliptic
partial differential operator parameterized by the control variable $\alpha
\in\Gamma\subset\mathbf{R}^{d}$,
\[
\mathcal{L}\left(  x,\alpha\right)  =\frac{1}{2}\sum_{\mu,\nu=1}^{n}\left(
\sum_{j=1}^{d}\sigma_{\mu j}\sigma_{\nu j}\left(  x\right)  \alpha_{j}%
^{2}\right)  \cdot\partial_{x^{\mu}x^{\nu}}+\sum_{i=1}^{n}\left(  b_{i}\left(
x\right)  +\sum_{j=1}^{d}h_{ij}(x)\alpha_{j}^{2}\right)  \partial_{x^{i}}%
+\sum_{j=1}^{d}f_{j}\left(  x,u\right)  \alpha_{j}^{2}.
\]

The fact that BSDEs on the space of linear probability could provide
probabilistic solutions for second order qusi-linear partial differential
equation (PDE) has been studied in \textrm{\cite{pp1,pp2,p1, p2, tang, yzhou,
ma}}. This probabilistic method was also extended to stochastic PDEs
\textrm{\cite{p10,pp2, pw}}, fully nonlinear cases \textrm{\cite{p3,ch, stz}}.
Recently, Peng \textrm{\cite{p5,p6,p7,p8,p9}} proposed the notion of sublinear
G-expectation and established associated stochastic calculus. In fact super
linear G-expectation can be introduced similarly. BSDE under super linear
G-expectation (G-BSDE) is also well defined under Lipschitz condition.
Generally, a new kind of BSDE corresponds to a new kind of PDE. Then a natural
question comes up: for what kind of PDE, BSDE under super linear G-expectation
can provide a probabilistic solution. This subject is also called Feynman-Kac
formula \textrm{\cite{kac,p2}}. Initially, we are not sure whether G-BSDE
corresponds to the HJB equation. This is the first reason that we write the
present paper. Secondly, we show that G-BSDE indeed provides a new
probabilistic solution of a class of HJB equations (see Peng \textrm{\cite{p3}
}for another interpretation). This fact shows that probabilistic
interpretations for solutions of HJB equation are not unique. Thirdly and more
importantly, the super (sub) linear G-expectation is in fact a super (sub)
strategy with mutually singular probability measures on the set of possible
paths, i.e., $\mathbb{E}^{\ast}\left[  \cdot\right]  =\underset{P\in
\mathcal{P}}{\sup}{\normalsize E}_{P}\left[  \cdot\right]  $, $\mathbb{E}%
_{\ast}\left[  \cdot\right]  =\underset{P\in\mathcal{P}}{\inf}{\normalsize E}%
_{P}\left[  \cdot\right]  $, where $\mathcal{P}$ is the set of risk-neutral
probabilities. An already known convenient framework to deal with super
strategy is the stochastic control framework. However, the connection between
the super strategy problem and stochastic control is not that obvious. Recall
that the stochastic control problem is the maximization of an expectation over
a set of processes. By our Feynman-Kac formula, we characterize the
G-expectation (G-BSDE) as a class of HJB equations, which establishes a direct
and equivalent relation among super strategy (G-expectation), HJB equation and
value function of a stochastic control. Superlinear G-expectation itself is
important in the theory of nonlinear expectation. We develop several
propositions which give deeper insight into properties of G-normal
distribution and G-Brownian motion, such as Proposition 2.2$\backsim$2.5.
Especially in Proposition 2.4, it is proved that the quadratic variance
$(\left\langle B\right\rangle _{t})$ of G-Brownian motion $(B_{t})$ is
differential in the sense of \textquotedblleft quasi-surely\textquotedblright%
\ for each $t$. The dynamic programming principle is also easily obtained in
our framework.

Motivated by the bid price in uncertainty volatility model \textrm{\cite{alp}%
}, this paper first studies super linear G-expectation. Sublinear
G-expectation and It\^{o} calculus under which have been well studied in Peng
\textrm{\cite{p5,p6,p7,p8,p9}}. We naturally want to know whether there is new
calculus under super linear G-expectation $\mathbb{E}_{\ast}$. However we show
that $\mathbb{E}^{\ast}$-Brownian motion is also a Brownian motion under
$\mathbb{E}_{\ast}$ and every super linear expectation $\mathbb{E}_{\ast}$ is
not dominated by itself but by sublinear expectation $\mathbb{E}^{\ast}%
[\cdot]:=-\mathbb{E}_{\ast}[-\cdot\ ]$. Due to the non-dominated property of
$\mathbb{E}_{\ast}$, we have to work with associated sublinear G-expectation
$\mathbb{E}^{\ast}$. In fact superlinear expectation and sublinear expectation
are complementary to each other. Super linear G-expectation is an auxiliary
means of sublinear G-expectation. It provides more insights to G-Brownian
motion and the uncertainty of random variables under sublinear expectation.
For instance, a `symmetric martingale' $(M_{t})$ \textrm{\cite{xz} }can be
characterized as that $(M_{t})$ is a martingale both under $\mathbb{E}^{\ast}$
and $\mathbb{E}_{\ast}$. Recently Li and Peng \textrm{\cite{lp}} established a
new framework for It\^{o} integral and related stochastic calculus. However,
there are many important and interesting problems still holding open under
this new framework, e.g., under what condition $\int_{0}^{\cdot}\eta_{s}%
dB_{s}$ is a martingale or a local martingale. So we will work within Peng's
framework, which are enough for us to illustrate our subject.

We first recall some notions under nonlinear expectation and proved that a
random variable $X$ is $\mathbb{E}_{\ast}$-normal distributed if and only if
$u(t,x)=\mathbb{E}_{\ast}\left[  \phi(x+\sqrt{t}X)\right]  $ is the viscosity
solution of PDE $\partial_{t}u-\frac{1}{2}\underset{\gamma\in\Gamma}{\inf
}tr\left\{  \gamma\gamma^{T}D^{2}u\right\}  =0$, $u\left(  0,x\right)
=\phi(x)$. Some useful property related to super linear G-expectation are also
listed. Section 3 establishes a uniqueness and existence theorem and a
comparison for BSDEs under super linear G-expectation. Section 4 studies the
relation between BSDEs under $\mathbb{E}_{\ast}$ and associated HJB equation.
In section 5 we show some applications: BSDEs under $\mathbb{E}_{\ast}$ could
characterize a class of stochastic control problems; we also give
representation for solutions of Black-Scholes-Barrenblett equations which
arise in mathematical finance and G-heat equations which are fundamentally
important in the theory of G-expectation. In the appendix we discuss the
dominated convergence theorem under sublinear expectation induced by mutually
singular probability measures.

\section{Super linear expectation: another point of view of sublinear
expectation}

For a given positive integer $n$ we will denote by $\left\langle
x,y\right\rangle $ the scalar product of $x,y\in\mathbf{R}^{n}$ and by
$\left\vert x\right\vert =\left\langle x,x\right\rangle ^{1/2}$ the Euclidean
norm of $x$. For two stochastic processes $(X_{t})$ and $(Y_{t})$, let
$\left\langle X,Y\right\rangle _{t}$ denote their mutual variance. We denote
by $\mathbb{S}(n)$ the collection of $n\times n$ symmetric matrices. We
observe that $\mathbb{S}(n)$ is an Euclidean space with the scalar product
$\left\langle A,B\right\rangle =tr[AB]$. Let $\Omega$ be a given set and let
${\mathcal{H}}$ be a linear space of real functions defined on $\Omega$ such
that if $X_{1},\ldots,X_{n}\in{\mathcal{H}}$ then $\varphi(X_{1},\ldots
X_{n})\in{\mathcal{H}}$ for each $\varphi\in$C$_{l.Lip(\mathbf{R}^{n})}$ where
$C_{l.Lip(\mathbf{R}^{n})}$ denotes the linear space of (local Lipschitz)
functions $\varphi$ satisfying
\[
|\varphi\left(  x\right)  -\varphi\left(  y\right)  |\leq C(1+\left\vert
x\right\vert ^{m}+\left\vert y\right\vert ^{m})|x-y|,\ \forall x,y\in
\mathbf{R}^{n},
\]
for some $C>0$, $m\in N$ depending on $\varphi$. ${\mathcal{H}}$ is considered
as a space of `random variables'. In this case $X=(X_{1},\ldots,X_{n})$ is
called an $n$-dimensional random vector, denoted by $X\in\mathcal{H}^{n}$. We
also denote by $\mathcal{B}(\Omega)$ the Borel $\sigma$-algebra of $\Omega$;
$C_{b}^{k}(\mathbf{R}^{n})$ the space of bounded and $k$-time continuously
differentiable functions with bounded derivatives of all orders less than or
equal to $k$; $C_{Lip(\mathbf{R}^{n})}$ the space of Lipschitz continuous functions.

\begin{definition}
A \textbf{nonlinear expectation} $\mathbb{E}$ on ${\mathcal{H}}$ is a
functional $\mathbb{E}:{\mathcal{H\mapsto}}\mathbf{R}$ satisfying the
following properties: for all $X,Y\in{\mathcal{H}}$, we have

(a) Monotonicity: If $X\geq Y$, then $\mathbb{E[}X\mathbb{]}\geq
\mathbb{E[}Y\mathbb{]}$.

(b) Constant preserving: $\mathbb{E[}c\mathbb{]=}c$.

If a functional $\mathbb{E}^{\ast}:{\mathcal{H\mapsto}}\mathbf{R}$ satisfies
(a), (b) and the following

(c) Sub-additivity: $\mathbb{E[}X+Y\mathbb{]}\leq\mathbb{E[}X\mathbb{]}%
+\mathbb{E[}Y\mathbb{]}$.

(d) Positive homogeneity: $\mathbb{E[\lambda}X\mathbb{]=\lambda E[}%
X\mathbb{]},\forall\lambda\geq0$.

then we call $\mathbb{E}^{\ast}$ a \textbf{sublinear expectation}. If a
functional $\mathbb{E}_{\ast}:{\mathcal{H\mapsto}}\mathbf{R}$ satisfies (a),
(b), (d) and

(c') Super-additivity: $\mathbb{E[}X+Y\mathbb{]}\geq\mathbb{E[}X\mathbb{]}%
+\mathbb{E[}Y\mathbb{]}$.

we call $\mathbb{E}_{\ast}$ a \textbf{superlinear expectation}.
\end{definition}

\begin{definition}
Let $X_{1}$ and $X_{2}$ be two n-dimensional random vectors defined on
nonlinear expectation spaces $(\Omega_{1},\mathcal{H}_{1},\mathbb{E}_{1})$ and
$( \Omega_{2},\mathcal{H}_{2},\mathbb{E}_{2})$ respectively. They are called
identically distributed, denoted by $X_{1}\overset{d}{=}X_{2}$, if%
\[
\mathbb{E}_{1}[\varphi(X_{1})]=\mathbb{E}_{2}[\varphi(X_{2})],\ \forall
\varphi\in C_{l.Lip}(\mathbf{R}^{n}).
\]

\end{definition}

\begin{definition}
In a nonlinear expectation space $(\Omega,\mathcal{H},\mathbb{E})$ a random
vector $Y\in{\mathcal{H}}^{n}$ is said to be independent of another random
vector $X\in{\mathcal{H}}^{m}$ under $\mathbb{E}$ if for each test function
$\varphi\in C_{l.Lip}(\mathbf{R}^{m+n})$ we have%
\[
\mathbb{E}[\varphi(X,Y)]=\mathbb{E}\left[  \mathbb{E}[\varphi(x,Y)]_{x=X}%
\right]  .
\]

\end{definition}

\begin{remark}
If $Y$ is independent of $X$, then $\mathbb{E[}X+Y\mathbb{]}=\mathbb{E[}%
X\mathbb{]}+\mathbb{E[}Y\mathbb{]}$.
\end{remark}

\begin{definition}
($G$-normal distribution with zero mean under positively homogeneous
expectation). A d-dimensional random vector $X=(X_{1},...,X_{d})$ in a
positively homogeneous expectation space $(\Omega,\mathcal{H},\mathbb{E})$ is
called G-normal distributed if for each $a,b\geq0$ we have%
\[
aX+b\overset{\_}{X}\,\overset{d}{=}\sqrt{a^{2}+b^{2}}X,
\]
where $\overset{\_}{X}$ is an independent copy of X.
\end{definition}

\begin{remark}
It is easy to check that $\mathbb{E}[X]=$ $\mathbb{E}[-X]=0$. The so called
`G' is related to $G:\mathbb{S}(d)\mapsto\mathbf{R}$ defined by
\[
G\left(  A\right)  =\frac{1}{2}\mathbb{E[}\left\langle AX,X\right\rangle
\mathbb{]},
\]

\end{remark}

\begin{proposition}
Let $\mathbb{E}_{\ast}$ be a super linear expectation. $\mathbb{E}^{\ast
}[\cdot]:=-\mathbb{E}_{\ast}[-\cdot\ ]$. Then

(i) $\mathbb{E}^{\ast}[\cdot]$ is a sublinear expectation.

(ii) There exist a family of linear expectation \{$E_{P};P\in\mathcal{P}$\} on
$(\Omega,\mathcal{H})$ such that $\mathbb{E}^{\ast}\left[  \cdot\right]
=\underset{P\in\mathcal{P}}{\sup}{\normalsize E}_{P}\left[  \cdot\right]  $,
$\mathbb{E}_{\ast}\left[  \cdot\right]  =\underset{P\in\mathcal{P}}{\inf
}{\normalsize E}_{P}\left[  \cdot\right]  $.

(iii) $G_{\ast}\left(  A\right)  :=\frac{1}{2}\mathbb{E_{\ast}[}\left\langle
AX,X\right\rangle \mathbb{]=}\frac{1}{2}\underset{\gamma\in\Gamma}{\inf
}tr\left\{  \gamma\gamma^{T}A\right\}  $, where $\Gamma$ is a bounded and
closed subset of $\mathbf{R}^{d\times d}$.

(iv) $\left\vert \mathbb{E}_{\ast}[X|\mathcal{F}_{t}]-\mathbb{E}_{\ast
}[Y|\mathcal{F}_{t}]\right\vert \leq\mathbb{E}^{\ast}[\left\vert
X-Y\right\vert |\mathcal{F}_{t}]$.
\end{proposition}

\paragraph{Proof}

It is easy to check (i). By properties of sublinear expectation
(\textrm{\cite{p5,hp}}), we obtain (ii) and (iii). (iv) is from Peng
\textrm{\cite{p4}.}

\begin{definition}
We introduce the natural Choquet capacity
\[
C^{\ast}(A):=\underset{P\in\mathcal{P}}{\sup}P(A),\ A\in\mathcal{B}(\Omega).
\]
A property holds \textquotedblleft quasi-surely\textquotedblright\ (q.s.) if
it holds outside a polar set A, i.e., $C^{\ast}(A)=0$. A mapping $X$ on
$\Omega$ with values in a topological space is said to be quasi-continuous
(q.c.) if $\forall\varepsilon>0$, there exists an open set $O$ with $C^{\ast
}(O)<\varepsilon$ such that $X|_{O^{c}}$ is continuous.
\end{definition}

\begin{definition}
(Viscosity solution). $u\in C([0,\infty)\times\mathbf{R}^{n})$ is called a
viscosity subsolution (supersolution) of backward PDE%
\begin{equation}
\partial_{t}u+F\left(  t,x,u,Du,D^{2}u\right)  =0,\ u\left(  T,x\right)
=\Phi(x).
\end{equation}
if for any $\varphi\in C_{b}^{1,3}([0,\infty]\times\mathbf{R}^{n})$, any
$(t,x)\in\lbrack0,\infty)\times\mathbf{R}^{n}$ which is a point of local
maximum of $u-\varphi$,
\begin{equation}
\partial_{t}\varphi+F\left(  t,x,\varphi,D\varphi,D^{2}\varphi\right)
\geq0,\ (resp.\leq0).
\end{equation}
$u$ is called a viscosity solution of PDE (2.1) if it is both super and subsolution.
\end{definition}

\begin{remark}
If PDE (2.1) is in forward form \textrm{\cite{cil}}, the signs `$+$', `$\geq
$', `$\leq$' in (2.2) are changed into `$-$', `$\leq$', `$\geq$' respectively.
Compared with definition of viscosity solution in \textrm{\cite{cil,ish,p3}},
we replace $C^{1,2}$ by $C_{b}^{1,3}$ just for technical convenience. In fact
there are trivial difference between these two definitions.
\end{remark}

\begin{proposition}
A random variable $X$ is $G_{\ast}$-normal distributed if and only if
\begin{equation}
u(t,x)=\mathbb{E}_{\ast}\left[  \varphi(x+\sqrt{t}X)\right]  ,(t,x)\in
\lbrack0,\infty)\times\mathbf{R}^{d},\varphi\in C_{l.Lip}(\mathbf{R}^{n}),
\end{equation}
is the viscosity solution of PDE
\begin{equation}
\partial_{t}u-G_{\ast}\left(  D^{2}u\right)  =0,u\left(  0,x\right)
=\varphi(x).
\end{equation}
where $D^{2}u=(\partial x^{i}x^{j}u)_{i,j=1}^{d}$, $G_{\ast}\left(  A\right)
=\frac{1}{2}\mathbb{E_{\ast}[}\left\langle AX,X\right\rangle \mathbb{]=}$
$\frac{1}{2}\underset{\gamma\in\Gamma}{\inf}tr\left\{  \gamma\gamma
^{T}A\right\}  $, $\Gamma$ is a bounded and closed subset of $\mathbf{R}%
^{d\times d}$.
\end{proposition}

\paragraph{Proof}

\bigskip We replace $t$ by $T-t$ in (2.3) and PDE(2.4), then the necessity is
an immediate result of Theorem 4.1 in the present paper. It can be also proved
similarly as in Peng \textrm{\cite{p6,p9}}.

Conversely, let $Y$ be a $G_{\ast}$-normal distributed random variable. Then
$\forall\varphi\in C_{l.Lip}(\mathbf{R}^{n})$, $v(t,x)=\mathbb{E}_{\ast
}\left[  \varphi(x+\sqrt{t}Y)\right]  $ is a viscosity solution of PDE (2.4).
If $\forall\varphi\in C_{l.Lip}(\mathbf{R}^{n})$, $u(t,x)=\mathbb{E}_{\ast
}\left[  \varphi(x+\sqrt{t}X)\right]  $ is also a viscosity solution of PDE
(2.4), then by the uniqueness of viscosity solution of PDE (2.4) (see
\textrm{\cite{cil,ish}}), we get that $v(t,x)=u(t,x)$, $\forall(t,x)\in
\lbrack0,\infty)\times\mathbf{R}^{d}$. Therefore we obtain that $X\,\overset
{d}{=}Y$, thus $X$ is $G_{\ast}$-normal distributed. $\Box$

\begin{remark}
The above proposition also holds under the framework of sublinear expectation.
\end{remark}

\begin{definition}
(\textbf{G--Brownian motion}). A d-dimensional process $(B_{t})_{t\geq0}$ on a
nonlinear expectation space $(\Omega,\mathcal{H},\mathbb{E})$ is called a
G--Brownian motion if the following properties are satisfied:
\end{definition}

(i) $B_{0}(\omega)=0$;

(ii) For each $t,s\geq0$, the increment $B_{t+s}-B_{t}$ is independent from
$(B_{t_{1}},B_{t_{2}},\ldots,B_{t_{n}})$, for each $n\in N$ and $0\leq
t_{1}\leq\cdots\leq t_{n}\leq t$;

(iii) $B_{t+s}-B_{t}\overset{d}{=}\sqrt{s}X$, where $X$ is G-normal distributed.

\begin{proposition}
Let $(B_{t})$ be a one-dimensional $G^{\ast}-$Brownian motion. $(\left\langle
B\right\rangle _{t})$ denotes its quadratic variance. Then

(i) $(\left\langle B\right\rangle _{t})$ is a continuous, increasing process
with finite variance, independent and stationary increments under
$\mathbb{E}^{\ast}$.

(ii)
\begin{equation}
\mathbb{E}^{\ast}[\varphi(\left\langle B\right\rangle _{t+s}-\left\langle
B\right\rangle _{s})|\mathcal{F}_{s}]=\underset{\underline{\sigma}^{2}%
\leq\alpha^{2}\leq\overline{\sigma}^{2}}{\sup}\varphi(\alpha^{2}%
t),\forall\varphi\in C(\mathbf{R}).
\end{equation}%
\begin{equation}
\mathbb{E}_{\ast}[\varphi(\left\langle B\right\rangle _{t+s}-\left\langle
B\right\rangle _{s})|\mathcal{F}_{s}]=\underset{\underline{\sigma}^{2}%
\leq\alpha^{2}\leq\overline{\sigma}^{2}}{\inf}\varphi(\alpha^{2}%
t),\forall\varphi\in C(\mathbf{R}).
\end{equation}
where wee denote the usual parameters $\overline{\sigma}^{2}=\mathbb{E}^{\ast
}[\left\langle B\right\rangle _{t=1}]$, $\underline{\sigma}^{2}=\mathbb{E}%
_{\ast}[\left\langle B\right\rangle _{t=1}]$.

(iii) For each $0\leq t\leq T<\infty$, we have
\begin{equation}
\underline{\sigma}^{2}(T-t)\leq\left\langle B\right\rangle _{T}-\left\langle
B\right\rangle _{t}\leq\overline{\sigma}^{2}(T-t)q.s..
\end{equation}
Therefore for $(\eta_{s})$ such that $\int_{0}^{t}\left\vert \eta
_{s}\right\vert ^{2}ds<\infty,\ q.s.$, we have
\begin{equation}
\underline{\sigma}^{2}\int_{0}^{t}\left\vert \eta_{s}\right\vert ^{2}%
ds\leq\int_{0}^{t}\left\vert \eta_{s}\right\vert ^{2}d\left\langle
B\right\rangle _{s}\leq\overline{\sigma}^{2}\int_{0}^{t}\left\vert \eta
_{s}\right\vert ^{2}ds,\ q.s..
\end{equation}

\end{proposition}

\paragraph{Proof}

See Peng \textrm{\cite{p5,p6,p8,p9}} for (i). $\forall\varphi\in
C_{l.Lip}(\mathbf{R})$, (2.5) and (2.6) hold true. Note that any $\varphi\in
C(\mathbf{R})$ can be approximated by $\varphi_{n}\in C_{l.Lip}(\mathbf{R})$
uniformly on a bounded subset of $\mathbf{R.}$ Thus $\forall\varphi\in
C(\mathbf{R})$, we have $\mathbb{E}_{\ast}[\varphi(\left\langle B\right\rangle
_{t+s}-\left\langle B\right\rangle _{s})|\mathcal{F}_{s}]=\underset
{n\rightarrow\infty}{\lim}\mathbb{E}_{\ast}[\varphi_{n}(\left\langle
B\right\rangle _{t+s}-\left\langle B\right\rangle _{s})|\mathcal{F}%
_{s}]=\underset{n\rightarrow\infty}{\lim}\underset{\underline{\sigma}^{2}%
\leq\alpha^{2}\leq\overline{\sigma}^{2}}{\min}\varphi_{n}(\alpha
^{2}t)=\underset{\underline{\sigma}^{2}\leq\alpha^{2}\leq\overline{\sigma}%
^{2}}{\min}\varphi(\alpha^{2}t)$. (2.7) is from \textrm{\cite{p9}}. (2.8) is a
sequence of (2.7). $\Box$

\begin{proposition}
Let $\mathbb{E}_{\ast}$ be a super linear expectation. $\mathbb{E}^{\ast
}[\cdot]:=-\mathbb{E}_{\ast}[-\cdot\ ]$. Let $(B_{t})$ be a one-dimensional
$G^{\ast}-$Brownian motion under $\mathbb{E}^{\ast}$. Then

(i) $(B_{t})$ and $(-B_{t})$ are Brownian motion under $\mathbb{E}_{\ast}$.

(ii) For each linear expectation $E_{P},P\in\mathcal{P}$ on $(\Omega
,\mathcal{H})$, $(B_{t})$ is a $E_{P}-$martingale and there exist an
$\mathcal{F}_{t}^{W}-$adapted process $(z_{t})$ such that $E_{P}\int_{0}%
^{t}\left\vert z_{s}\right\vert ^{2}ds<\infty$ and $B_{t}=\int_{0}^{t}%
z_{s}dW_{s}$, $(W_{t})$ is a standard $E_{P}-$Brownian motion with
$\underline{\sigma}^{2}\leq z_{t}^{2}\leq\overset{\_}{\sigma}^{2}$, for a.e.
t, $P-a.s..$

(iii) $\frac{d\left\langle B\right\rangle _{t}}{dt}$ exists q.s. for each
$t\geq0$ and $\underline{\sigma}^{2}\leq\frac{d\left\langle B\right\rangle
_{t}}{dt}\leq\overset{\_}{\sigma}^{2}$.

(iv) $\mathbb{E}^{\ast}[\varphi(\frac{d\left\langle B\right\rangle _{t}}%
{dt})]=\underset{P}{\sup}\varphi(\upsilon_{P})$, $\forall\varphi\in
C(\mathbf{R})$.
\end{proposition}

\paragraph{Proof}

\bigskip(i) We only check (iii) in Definition 2.7. Since $(B_{t+s}%
-B_{t})\overset{d}{=}\sqrt{s}X$, $X$ is a $G^{\ast}-$normal distributed random
variable. So $\mathbb{E}_{\ast}[\varphi(B_{t+s}-B_{t})]=-\mathbb{E}^{\ast
}[-\varphi(B_{t+s}-B_{t})]=-\mathbb{E}^{\ast}[-\varphi(\sqrt{s}X)]=\mathbb{E}%
_{\ast}[\varphi(\sqrt{s}X)]$, $\forall\varphi\in C_{l.Lip}(\mathbf{R})$. On
the other hand, it is easy to check that $X$ is also $\mathbb{E}_{\ast}%
-$normal distributed. By Definition 2.7, $(B_{t})$ is a $\mathbb{E}_{\ast}%
$-Brownian motion. Similarly $(-B_{t})$ is also a Brownian motion under
$\mathbb{E}_{\ast}$.

(ii) It is easy seen that $\mathbb{E}^{\ast}\left[  \cdot|\mathcal{F}%
_{t}\right]  =\underset{P\in\mathcal{P}}{\sup}{\normalsize E}_{P}\left[
\cdot|\mathcal{F}_{t}\right]  $ and $\mathbb{E}_{\ast}\left[  \cdot
|\mathcal{F}_{t}\right]  =\underset{P}{\inf}{\normalsize E}_{P}\left[
\cdot|\mathcal{F}_{t}\right]  $. For each $P\in\mathcal{P}$, $E_{P}%
[B_{t+s}-B_{t}|\mathcal{F}_{t}]\leq\mathbb{E}^{\ast}[B_{t+s}-B_{t}%
|\mathcal{F}_{t}]=0$ and $E_{P}[B_{t+s}-B_{t}|\mathcal{F}_{t}]\geq
\mathbb{E}_{\ast}[B_{t+s}-B_{t}|\mathcal{F}_{t}]=0$. So we have $E_{P}%
[B_{t+s}|\mathcal{F}_{t}]=B_{t}$. By the classical martingale representation,
there exist an $\mathcal{F}_{t}^{W}-$adapted process $(z_{t})$ such that
$E_{P}\int_{0}^{t}\left\vert z_{s}\right\vert ^{2}ds<\infty$ and $B_{t}%
=\int_{0}^{t}z_{s}dW_{s}$, $(W_{t})$ is a standard $E_{P}-$Brownian motion. By
Proposition 2.3(iii), we have $\underline{\sigma}^{2}t\leq\left\langle
B\right\rangle _{t}=\int_{0}^{t}\left\vert z_{s}\right\vert ^{2}ds\leq$
$\overline{\sigma}^{2}t$ for any $t$, thus $\underline{\sigma}^{2}%
\leq\left\vert z_{t}\right\vert ^{2}\leq$ $\overline{\sigma}^{2}$.

(iii) By (ii) we know that for each $P\in\mathcal{P}$, $\left\langle
B\right\rangle _{t}$ is $P$-a.s. differential at any $t$. For each given $t$,
let $A$ denotes the set of $\omega\in\Omega$ such that$\underset{s\rightarrow
t}{\ \lim}\frac{\left\langle B^{j}\right\rangle _{s}-\left\langle
B^{j}\right\rangle _{t}}{s-t}$ does not exist. Then we obtain that for each
$P\in\mathcal{P}$, $P(A)=0$ and $C^{\ast}(A)=\underset{P\in\mathcal{P}}{\sup
}P(A)=0$. Therefore $\frac{d\left\langle B\right\rangle _{t}}{dt}$ exists
\textit{q.s.} for each $t\geq0$.

(iv) By (ii), for each $P\in\mathcal{P}$, $\ $there is a constant
$\upsilon_{P}\in\lbrack\underline{\sigma}^{2},\overset{\_}{\sigma}^{2}]$ such
that $E_{P}[\varphi(\frac{d\left\langle B\right\rangle _{t}}{dt}%
)]=E_{P}[\varphi(z_{t}^{2})]=\varphi(\upsilon_{P})$,$\ $so $\mathbb{E}^{\ast
}[\varphi(\frac{d\left\langle B\right\rangle _{t}}{dt})]=\underset
{P\in\mathcal{P}}{\sup}{\normalsize E}_{P}\left[  \varphi(\frac{d\left\langle
B\right\rangle _{t}}{dt})\right]  =\underset{P}{\sup}\varphi(\upsilon_{P})$,
$\forall\varphi\in C(\mathbf{R})$. $\Box$

Every G-Brownian motion is related to a sub (super) linear function G and a
nonempty, bounded and closed subset $\Gamma$ of $\mathbf{R}^{d\times d}$.
$\Gamma$ characterizes the variance uncertainty of corresponding G-Brownian
motion. One should note that for a general $\Gamma$, components of $G_{\Gamma
}$-Brownian motion may be not independent from others. So the item
$\left\langle B^{i},B^{j}\right\rangle _{t}$ arises in the G-It\^{o} formula
and G-stochastic differential equation.

Consider the following typical nonlinear heat equation:%

\begin{equation}
\partial_{t}u(t,x)-\frac{1}{2}\sum_{i=1}^{d}\left[  \underline{\sigma}_{i}%
^{2}(\partial_{x^{i}x^{i}}u(t,x))^{+}-\overline{\sigma}_{i}^{2}(\partial
_{x^{i}x^{i}}u(t,x))^{-}\right]  =0.
\end{equation}
The above equation corresponds to $\Gamma=\{diag[\gamma_{1},\ldots,\gamma
_{d}]$, $\gamma_{i}^{2}\in\lbrack\underline{\sigma}_{i}^{2},\overline{\sigma
}_{i}^{2}],i=1,\ldots d\}$. Let $(B_{t})=(B_{t}^{1},\ldots,B_{t}^{d})$ be a
$d$-dimensional G-Brownian motion, then

\begin{proposition}
If $(B_{t}^{j})$ is independent of $(B_{t}^{i})$, $1\leq i<j\leq d$, then

(i) $u(t,x)=\mathbb{E}_{\ast}\left[  \phi(x+B_{t})\right]  $ is the viscosity
solution of PDE (2.8) with $u\left(  0,x\right)  =\phi(x)$.

(ii) $\Gamma=\{diag[\gamma_{1},\ldots,\gamma_{d}],\gamma_{i}^{2}\in
\lbrack\underline{\sigma}_{i}^{2},\overline{\sigma}_{i}^{2}],i=1,\ldots d\}$.
\end{proposition}

\paragraph{Proof}

\bigskip For convenience, we only consider two dimensional G-Brownian motion.
If $(B_{t}^{2})$ is independent of $(B_{t}^{1})$, then similarly to
Proposition 2.2, we can prove that $u(t,x_{1},x_{2})=\mathbb{E}_{\ast}\left[
\phi(x_{1}+B_{t}^{1},x_{2}+B_{t}^{2})\right]  $ is the viscosity solution of
PDE (2.8), which leads to that $\Gamma$ is a set of diagonal matrices. $\Box$

\bigskip Let $\Omega=C_{0}^{d}(\mathbf{R}^{+})$ denote the space of all
$\mathbf{R}^{d}-$valued continuous paths $(\omega_{t})_{t\in R^{+}}$ with
$\omega_{0}=0$, by $C_{b}(\Omega)$ all bounded and continuous functions on
$\Omega$. For each fixed $T\geq0$, we consider the following space of random
variables:%
\[
C_{l.Lip}(\Omega_{T}):=\{X(\omega)=\varphi(\omega_{t_{1}\wedge T}%
,...,\omega_{t_{m}\wedge T}),\forall m\geq1,\forall\varphi\in
C_{l.Lip(\mathbb{R}^{m})}\}.
\]
We also denote
\[
C_{l.Lip}(\Omega):=\overset{\infty}{\underset{n=1}{\cup}}C_{l.Lip}(\Omega
_{n}).
\]
We will consider the canonical space and set $B_{t}(\omega)=\omega_{t}$. For a
given sublinear function $G^{\ast}\left(  A\right)  =$ $\frac{1}{2}%
\underset{\gamma\in\Gamma}{\sup}tr\left\{  \gamma\gamma^{T}A\right\}  $, where
$A\in$ $\mathbb{S}(d)$, $\Gamma$ is a given nonempty, bounded and closed
subset of $\mathbf{R}^{d\times d}$, by the following%
\[
\partial_{t}u(t,x)-G^{\ast}\left(  D_{x}^{2}u\right)  =0,\ u(0,x)=\varphi(x),
\]
Peng \textrm{\cite{p5} }defined\textrm{ }$G^{\ast}$-expectation $\mathbb{E}%
^{\ast}$ as $\mathbb{E}^{\ast}[\varphi(x+B_{t})]=u(t,x)$. For each $p\geq1$,
$X\in C_{l.Lip}(\Omega)$, $\Vert X\Vert_{p}=\mathbb{E}^{\ast}[\left\vert
X\right\vert ^{p}]^{\frac{1}{p}}$ forms a norm and $\mathbb{E}^{\ast}$ can be
continuously extended to a Banach space, denoted by $L{_{G}^{\mathrm{p}}%
}(\Omega)$. Hu and Peng \textrm{\cite{hp} }proved that $L{_{G}^{\mathrm{p}}%
}(\Omega)=\{X|\ X$ is $C^{\ast}-$quasi-continuous, $\mathcal{B}(\Omega
)-$measurable function s.t. $\underset{P\in\mathcal{P}}{\sup}{\normalsize E}%
_{P}[\left\vert X\right\vert ^{p}]<\infty\}$. By the method of Markov chains,
Peng \textrm{\cite{p5,p7}} also defined corresponding conditional expectation,
$\mathbb{E}^{\ast}\left[  \cdot|{\mathcal{F}}_{t}\right]  :L{_{G}^{\mathrm{1}%
}}(\Omega)\mapsto L{_{G}^{\mathrm{1}}}(\Omega_{t})$, where ${\mathcal{F}}%
_{t}:=\mathcal{B}(\Omega_{t})$, $\Omega_{t}:=\left\{  \omega._{\wedge
t}:\omega\in\Omega\right\}  $. Under $\mathbb{E}^{\ast}\left[  \cdot\right]
$, the canonical process $B_{t}(\omega)=\omega_{t}$, $t\in\lbrack0,\infty)$ is
a $G^{\ast}$-Brownian motion. The $G^{\ast}$-expectation $\mathbb{E}^{\ast}$
can be extended to more general space. We define the upper expectation for
each $\mathcal{B}(\Omega)$-measure random variable which makes the following
definition meaningful:%
\[
\overset{\_}{\mathbb{E}}\left[  X\right]  =\underset{P\in\mathcal{P}}{\sup
}{\normalsize E}_{P}\left[  X\right]  .
\]
Note that $\overset{\_}{\mathbb{E}}=\mathbb{E}^{\ast}$ on $L{_{G}^{\mathrm{p}%
}}(\Omega)$ and $\forall A\in\mathcal{B}(\Omega)$, $\overset{\_}{\mathbb{E}%
}\left[  \mathbf{1}_{A}\right]  =C^{\ast}(A)$. Without the loss of generality,
we still denote $\overset{\_}{\mathbb{E}}$ by $\mathbb{E}^{\ast}$.

Similarly, we can define superlinear expectation, $\mathbb{E}_{\ast}\left[
\cdot\right]  $ and $\mathbb{E}_{\ast}\left[  \cdot|{\mathcal{F}}_{t}\right]
$ by the following PDE%
\begin{equation}
\partial_{t}u(t,x)-G_{\ast}\left(  D_{x}^{2}u\right)  =0,\ u(0,x)=\varphi(x).
\end{equation}
where $G_{\ast}\left(  A\right)  =\frac{1}{2}\underset{\gamma\in\Gamma}{\inf
}tr\left\{  \gamma\gamma^{T}A\right\}  $. It is easy to check that
$\mathbb{E}_{\ast}\left[  \cdot\right]  =$ $-\mathbb{E}^{\ast}[-\cdot\ ]$ and
$\mathbb{E}_{\ast}\left[  \cdot|{\mathcal{F}}_{t}\right]  =$ $-\mathbb{E}%
^{\ast}[-\cdot\ |{\mathcal{F}}_{t}]$. The following property hold for
$\mathbb{E}_{\ast}\left[  \cdot|{\mathcal{F}}_{t}\right]  $ \textit{q.s.}.

\begin{proposition}
For $X,Y\in L{_{G}^{\mathrm{1}}}( \Omega)$, we have q.s.,

(i) $\mathbb{E}_{\ast}[\eta X|\mathcal{F}_{t}]=\eta^{+}\mathbb{E}_{\ast
}[X|\mathcal{F}_{t}]+\eta^{-}\mathbb{E}_{\ast}[-X|\mathcal{F}_{t}]$, for
bounded $\eta\in L{_{G}^{\mathrm{1}}}(\Omega_{t})$.

(ii) If $\mathbb{E}_{\ast}[X|\mathcal{F}_{t}]=-\mathbb{E}_{\ast}%
[-X|\mathcal{F}_{t}]$, for some $t$, then $\mathbb{E}_{\ast}\mathbb{[}%
X+Y|\mathcal{F}_{t}\mathbb{]}=\mathbb{E}_{\ast}\mathbb{[}X|\mathcal{F}%
_{t}\mathbb{]}+\mathbb{E}_{\ast}\mathbb{[}Y|\mathcal{F}_{t}\mathbb{]}$.

(iii) $\mathbb{E}_{\ast}\mathbb{[}X+\eta|\mathcal{F}_{t}\mathbb{]}%
=\mathbb{E}_{\ast}\mathbb{[}X|\mathcal{F}_{t}\mathbb{]}+\eta$, $\eta\in
L{_{G}^{\mathrm{1}}}(\Omega_{t})$.
\end{proposition}

For a partition of $[0,T]$: $0=t_{0}<t_{1}<\cdots<t_{N}=T$, we set

${\mathcal{M}_{G}^{\mathrm{p,0}}}(0,{\mathit{T}})$: the collection of
processes $\eta_{t}(\omega)=\sum_{j=0}^{N}\xi_{j}(\omega)\cdot1_{[t_{j}%
,t_{j+1}]}(t)$, where $\xi_{j}(\omega)\in L{_{G}^{\mathrm{p}}}(\Omega_{t_{j}%
}),j=0,1,...,N.$

$\mathcal{M}{_{G}^{\mathrm{P}}}(0,{\mathit{T}})$: the completion of
${\mathcal{M}_{G}^{\mathrm{p,0}}}(0,{\mathit{T}})$ under norm $||\eta
||=\left(  \mathbb{E}^{\ast}\left[  \int_{0}^{T}|\eta_{t}|^{p}dt\right]
\right)  ^{^{\frac{1}{p}}}$.

$\mathcal{H}{_{G}^{\mathrm{P}}}(0,{\mathit{T}})$: the completion of
${\mathcal{M}_{G}^{\mathrm{p,0}}}(0,{\mathit{T}})$ under norm $||\eta||_{\ast
}=\left(  \int_{0}^{T}\mathbb{E}^{\ast}[|\eta_{t}|^{p}]dt\right)  ^{^{\frac
{1}{p}}}$.

It\^{o} integral for process $\eta_{t}\in\mathcal{M}{_{G}^{\mathrm{2}}%
}(0,{\mathit{T}})$ is well defined in Peng \textrm{\cite{p9}}. A new framework
of It\^{o} integral is constructed in Li and Peng \textrm{\cite{lp}}. An
important difference between these two integral is that integrands in the
latter It\^{o} integral may be not quasi-continuous. It\^{o} formula have been
obtained in Peng \textrm{\cite{p5,p7,p9}}, Gao \textrm{\cite{gao}} and
generalized by Li and Peng \textrm{\cite{lp}. }We now adapt\textrm{ }Li and
Peng's formula to our framework. In fact Gao's and Peng's are enough for us to use.

\begin{proposition}
For an It\^{o} process $X_{t}=x_{0}+\int_{0}^{t}b_{s}ds+\int_{0}^{t}\eta
_{s}^{ij}d\left\langle B^{i},B^{j}\right\rangle _{s}+\int_{0}^{t}\beta_{s}%
^{j}dB_{s}^{j}$, where $x_{0}\in\mathbf{R}^{n}$, $b_{s},\eta_{s}^{i,j}%
,\beta_{s}^{j}\in\mathcal{M}{_{G}^{\mathrm{2}}}(0,{\mathit{T;}}\mathbf{R}%
^{n})$ and $\Phi\in C^{1,2}([0,\infty)\times\mathbf{R}^{n}\mapsto\mathbf{R})$
s.t $\partial\Phi_{t}$, $\partial_{x^{\nu}}\Phi b_{s}^{\nu}$, $\partial
_{x^{\nu}}\Phi\eta_{s}^{\nu ij}$, $\partial_{x^{\mu}x^{\nu}}\Phi\beta_{s}^{\nu
i}\beta_{s}^{\nu j}\in\mathcal{M}{_{G}^{\mathrm{1}}}(0,{\mathit{T;}}%
\mathbf{R})$, $\partial_{x^{\nu}}\Phi\beta_{s}^{\nu j}\in\mathcal{M}%
{_{G}^{\mathrm{2}}}(0,{\mathit{T;}}\mathbf{R})$, then for each $t$, we have
q.s.,%
\begin{align*}
\Phi\left(  t,X_{t}\right)   &  =\Phi\left(  0,x_{0}\right)  +\int_{0}%
^{t}\partial_{x^{\nu}}\Phi\beta_{s}^{\nu j}dB_{s}^{j}+\int_{0}^{t}%
(\partial\Phi_{t}+\partial_{x^{\nu}}\Phi b_{s}^{\nu})ds\\
&  +\int_{0}^{t}(\partial_{x^{\nu}}\Phi\eta_{s}^{\nu ij}+\frac{1}{2}%
\partial_{x^{\mu}x^{\nu}}\Phi\beta_{s}^{\nu i}\beta_{s}^{\nu j})d\left\langle
B^{i},B^{j}\right\rangle _{s}.
\end{align*}
Here we use the Einstein convention, i.e., the above repeated indices $\mu
,\nu,i,j$ imply summation.
\end{proposition}

\section{BSDEs under super linear expectation and comparison}

We consider the following $n$-dimensional backward stochastic differential
equation (BSDE)%

\begin{equation}
Y_{t}=\mathbb{E}_{\ast}\left[  \xi+\int_{t}^{T}g\left(  s,Y_{s}\right)
ds+\int_{t}^{T}f\left(  s,Y_{s}\right)  d\left\langle B\right\rangle
_{s}|{\mathcal{F}}_{t}\right]
\end{equation}
where

{\textbf{(H3.1)}}. $\xi\in L_{G}^{1}(\mathcal{F}_{T})$.

{\textbf{(H3.2)}}. There exists a constant $K\geq0$ , such that for $a.e.t$,
$q.s.$, $\forall y^{1},y^{2},z^{1},z^{2}$:$\quad$%
\[
|\mathit{g}\left(  t,y^{1}\right)  -\mathit{g}\left(  t,y^{2}\right)  |\leq
K|y^{1}-y^{2}|,\ |\mathit{f}\left(  t,y^{1}\right)  -\mathit{f}\left(
t,y^{2}\right)  |\leq K|y^{1}-y^{2}|,
\]

and the process $\left(  \mathit{g}\left(  t,\mathrm{0}\right)  \right)
_{t\in\left[  0,T\right]  }$ and $\left(  \mathit{f}\left(  t,\mathrm{0}%
\right)  \right)  _{t\in\left[  0,T\right]  }\in\mathcal{M}{_{G}^{\mathrm{1}}%
}(0,{\mathit{T}})$.

\begin{theorem}
\label{thm3.1} There is a unique solution $(Y_{t})\in{\mathcal{H}%
_{G}^{\mathrm{1}}}(0,{\mathit{T}};\mathbf{R}^{n})$ for (3.1).
\end{theorem}

\paragraph{Proof}

Consider a mapping $\Lambda$ from $\mathcal{M}{_{G}^{\mathrm{1}}}$ to
${\mathcal{H}_{G}^{\mathrm{1}}}$ defined as follows, for $Y\in\mathcal{M}%
{_{G}^{\mathrm{1}}}$,
\[
\Lambda_{t}(Y)=\mathbb{E}_{\ast}\left[  \xi+\int_{t}^{T}g\left(
s,Y_{s}\right)  ds+\int_{t}^{T}f\left(  s,Y_{s}\right)  d\left\langle
B\right\rangle _{s}|{\mathcal{F}}_{t}\right]  .
\]
It is easy to deduce from Lipschitz conditions that indeed $\Lambda_{t}%
(Y)\in{\mathcal{H}_{G}^{\mathrm{1}}}$. A solution of (3.1) is a fixed point of
the mapping $\Lambda.$ Uniqueness and existence of a fixed point will follow
the fact that , for each $T>0$, $\Lambda$ is a strict contraction on
${\mathcal{H}_{G}^{\mathrm{1,}\beta}}$ equipped with the norm $||X||_{\beta
}=\left[  \int_{0}^{T}e^{-2\beta t}\mathbb{E}^{\ast}|X_{t}|dt\right]  $ for
$\beta$ large enough.

Let $Y,Y^{\prime}\in{\mathcal{H}_{G}^{\mathrm{1}}}$, by Proposition 2.1(iv),
we have
\begin{align*}
\mathbb{E}^{\ast}|\Lambda_{t}(Y)-\Lambda_{t}(Y^{\prime})|  &  \leq
\mathbb{E}^{\ast}|\mathbb{E}_{\ast}\left[  \xi+\int_{t}^{T}g\left(
s,Y_{s}\right)  ds+\int_{t}^{T}f\left(  s,Y_{s}\right)  d\left\langle
B\right\rangle _{s}|{\mathcal{F}}_{t}\right] \\
&  -\mathbb{E}_{\ast}\left[  \xi+\int_{t}^{T}g\left(  s,Y_{s}^{\prime}\right)
ds+\int_{t}^{T}f\left(  s,Y_{s}^{\prime}\right)  d\left\langle B\right\rangle
_{s}|{\mathcal{F}}_{t}\right]  |\\
&  \leq\mathbb{E}^{\ast}|\int_{t}^{T}g\left(  s,Y_{s}\right)  -g\left(
s,Y_{s}^{\prime}\right)  ds+\int_{t}^{T}f\left(  s,Y_{s}\right)  -f\left(
s,Y_{s}^{\prime}\right)  d\left\langle B\right\rangle _{s}|\\
&  \leq\mathbb{E}^{\ast}\int_{t}^{T}|g\left(  s,Y_{s}\right)  -g\left(
s,Y_{s}^{\prime}\right)  |ds+\overline{\sigma}^{2}\mathbb{E}^{\ast}\int
_{t}^{T}|f\left(  s,Y_{s}\right)  -f\left(  s,Y_{s}^{\prime}\right)  |ds\\
&  \leq(1+\overline{\sigma}^{2})K\mathbb{E}^{\ast}\int_{t}^{T}|Y_{s}%
-Y_{s}^{\prime}|ds
\end{align*}

Let $\beta=K(1+\overline{\sigma}^{2})$. Multiply $e^{2\beta t}$ on both sides
of above inequality and then integrate them on $[0,T]$, we deduce that
$||\Lambda_{t}(Y)-\Lambda_{t}(Y^{\prime})||_{\beta}$ $\leq\frac{1}{2}%
||Y_{s}-Y_{s}^{\prime}||_{\beta}$. Note that ${\mathcal{H}_{G}^{\mathrm{1}}}$
is a sub space of $\mathcal{M}{_{G}^{\mathrm{1}}}$, thus by the contract
mapping principle and the equivalent of norms, there is a unique fixed point
for operator $\Lambda$ under norm $||\cdot||_{\ast}=\int_{0}^{T}%
\mathbb{E}^{\ast}|\cdot|dt$. $\Box$

\bigskip

For the following one dimensional linear BSDE:%

\begin{equation}
Y_{s}=\mathbb{E}_{\ast}\left[  \xi+\int_{t}^{T}\left(  a_{s}Y_{s}%
+A_{s}\right)  ds+\int_{t}^{T}\left(  b_{s}Y_{s}+C_{s}\right)  d\left\langle
B\right\rangle _{s}|\mathcal{F}_{t}\right]  ,\ t\in\lbrack0,T],
\end{equation}

Assuming {\textbf{(H3.1)}} and

{\textbf{(H3.2)}}$\prime$. There exists a constant $K\geq0$ , such that for
$a.e.t$, $q.s.$,$\quad$%
\[
|a_{t}|+|b_{t}|\leq K,
\]

and the process $\left(  A_{t}\right)  _{t\in\left[  0,T\right]  }$ and
$\left(  C_{t}\right)  _{t\in\left[  0,T\right]  }\in$ ${\mathcal{M}%
_{G}^{\mathrm{1}}}(0,{\mathit{T}};\mathbf{R})$.

\begin{theorem}
The following process
\begin{equation}
Y_{t}=Q_{t}^{-1}\mathbb{E}_{\ast}\left[  Q_{T}\xi+\int_{t}^{T}Q_{s}%
A_{s}ds+\int_{t}^{T}Q_{s}C_{s}d\left\langle B\right\rangle _{s}|\mathcal{F}%
_{t}\right]
\end{equation}
solves BSDE (3.2), where
\[
Q_{t}=\exp\{\int_{0}^{t}a_{s}ds+\int_{0}^{t}b_{s}d\left\langle B\right\rangle
_{s}\}.
\]

\end{theorem}

\paragraph{Proof}

By G-It\^{o} formula, we have%

\[
dQ_{t}=d\exp\{\int_{0}^{t}a_{s}ds+\int_{0}^{t}b_{s}d\left\langle
B\right\rangle _{s}\}=Q_{t}a_{t}dt+Q_{t}b_{t}d\left\langle B\right\rangle
_{t}.
\]

\bigskip We denote $M_{t}:=\mathbb{E}_{\ast}\left[  \xi+\int_{0}^{T}\left(
a_{s}Y_{s}+A_{s}\right)  ds+\int_{0}^{T}\left(  b_{s}Y_{s}+C_{s}\right)
d\left\langle B\right\rangle _{s}|\mathcal{F}_{t}\right]  $. Then by
2-dimensional G-It\^{o} formula (It holds similarly to the classical It\^{o}
formula \textit{w.r.t} semi-martingale) we get that%

\begin{align*}
d(Q_{t}Y_{t})  &  =Q_{t}dM_{t}-Q_{t}\left[  \left(  a_{t}Y_{t}+A_{t}\right)
dt+\left(  b_{t}Y_{t}+C_{t}\right)  d\left\langle B\right\rangle _{t}\right]
+Q_{t}a_{t}Y_{t}dt+Q_{t}b_{t}Y_{t}d\left\langle B\right\rangle _{t}\\
&  =Q_{t}dM_{t}-Q_{t}A_{t}dt-Q_{t}C_{t}d\left\langle B\right\rangle _{t}%
\end{align*}
Therefore%
\[
Q_{t}Y_{t}=Q_{0}Y_{0}+\int_{0}^{t}Q_{s}dM_{s}-\int_{0}^{t}Q_{s}A_{s}%
ds-\int_{0}^{t}Q_{s}C_{s}d\left\langle B\right\rangle _{s}%
\]
and
\[
Q_{T}Y_{T}=Q_{0}Y_{0}+\int_{0}^{T}Q_{s}dM_{s}-\int_{0}^{T}Q_{s}A_{s}%
ds-\int_{0}^{T}Q_{s}C_{s}d\left\langle B\right\rangle _{s}%
\]

Note that $\int_{0}^{t}Q_{s}dM_{s}$ is a $\mathbb{E}_{\ast}$-martingale,
therefore%
\begin{align*}
Q_{t}Y_{t}  &  =\mathbb{E}_{\ast}\left[  Q_{T}\xi+\int_{0}^{T}Q_{s}%
A_{s}ds+\int_{0}^{T}Q_{s}C_{s}d\left\langle B\right\rangle _{s}|{\mathcal{F}%
}_{t}\right]  -\int_{0}^{t}Q_{s}A_{s}ds-\int_{0}^{t}Q_{s}C_{s}d\left\langle
B\right\rangle _{s}\\
&  =\mathbb{E}_{\ast}\left[  Q_{T}\xi+\int_{t}^{T}Q_{s}A_{s}ds+\int_{t}%
^{T}Q_{s}C_{s}d\left\langle B\right\rangle _{s}|{\mathcal{F}}_{t}\right]
\end{align*}
and
\[
Y_{t}=Q_{t}^{-1}\mathbb{E}_{\ast}\left[  Q_{T}\xi+\int_{t}^{T}Q_{s}%
A_{s}ds+\int_{t}^{T}Q_{s}C_{s}d\left\langle B\right\rangle _{s}|{\mathcal{F}%
}_{t}\right]  .
\]
$\Box$

Let $(Y_{t})$, $(\overline{Y_{t}})$ be the unique solutions of the following
two one dimensional BSDEs:
\begin{equation}
Y_{t}=\mathbb{E}_{\ast}\left[  \xi+\int_{t}^{T}g\left(  s,Y_{s}\right)
ds+\int_{t}^{T}f\left(  s,Y_{s}\right)  d\left\langle B\right\rangle
_{s}|{\mathcal{F}}_{t}\right]
\end{equation}%
\begin{equation}
\overline{Y_{t}}=\mathbb{E}_{\ast}\left[  \overline{\xi}+\int_{t}^{T}%
\overline{g}\left(  s,\overline{Y_{s}}\right)  ds+\int_{t}^{T}\overline
{f}\ \left(  s,\overline{Y_{s}}\right)  d\left\langle B\right\rangle
_{s}|{\mathcal{F}}_{t}\right]
\end{equation}
We now establish a comparison between them.

\begin{theorem}
If $\xi\geq\overline{\xi}$, $q.s.$, and for $a.e.t$, $q.s.$, $\forall
y\in\mathbf{R}$, $\delta\mathbb{\geq}0$,%
\begin{equation}
g\left(  t,y+\delta\right)  \geq\overline{g}\left(  t,y\right)  \quad
\text{and}\quad f\left(  t,y+\delta\right)  \geq\overline{f}\left(
t,y\right)  ,
\end{equation}

then we have $\forall t\in[0,T]$, $Y_{t}\geq\overline{Y_{t}}$, q.s..
\end{theorem}

\paragraph{Proof}

Given any $y_{0}\in\mathbf{R}$, we define the following two sequences:%

\[
Y_{t}^{i+1}=\mathbb{E}_{\ast}\left[  \xi+\int_{t}^{T}g\left(  s,Y_{s}%
^{i}\right)  ds+\int_{t}^{T}f\left(  s,Y_{s}^{i}\right)  d\left\langle
B\right\rangle _{s}|{\mathcal{F}}_{t}\right]  ,
\]%
\[
\overline{Y_{t}}^{i+1}=\mathbb{E}_{\ast}\left[  \overline{\xi}+\int_{t}%
^{T}\overline{g}\left(  s,\overline{Y_{s}}^{i}\right)  ds+\int_{t}%
^{T}\overline{f}\ \left(  s,\overline{Y_{s}}^{i}\right)  d\left\langle
B\right\rangle _{s}|{\mathcal{F}}_{t}\right]  .
\]
Set $Y_{t}^{0}=\overline{Y_{t}}^{0}=y_{0}$. Obviously $\left\{  Y_{t}%
^{i}\right\}  _{i=0}^{\infty}$ and $\left\{  \overline{Y_{t}}^{i}\right\}
_{i=0}^{\infty}$ are Cauchy sequences in ${\mathcal{H}_{G}^{\mathrm{1}}%
}(0,{\mathit{T}};\mathbf{R})$ and $Y_{t}^{i}\rightarrow Y_{t}$, $\overline
{Y_{t}}^{i}\rightarrow\overline{Y_{t}}$.

One can check that, by the monotonicity of super linear expectation
$\mathbb{E}_{\ast}\left[  \cdot|{\mathcal{F}}_{t}\right]  $ , we have
$Y_{t}^{1}\geq$ $\overline{Y_{t}}^{1}$, $Y_{t}^{2}\geq$ $\overline{Y_{t}}^{2}%
$, $\cdots$, $Y_{t}^{i}\geq$ $\overline{Y_{t}}^{i}$, $\cdots$. Thus we get
that $Y_{t}=\lim\limits_{i\rightarrow\infty}Y_{t}^{i}\geq\lim
\limits_{i\rightarrow\infty}$ $\overline{Y_{t}}^{i}=\overline{Y_{t}}$. $\Box$

\section{Probabilistic interpretation for a class of HJB equations}

The Hamilton-Jacobi-Bellman (HJB) equation is a second order fully nonlinear
partial differential equation which is central to optimal control theory. The
solution of the HJB equation is the `value function', which gives the optimal
cost-to-go for a given dynamical system with an associated cost function. In
general case, the HJB equation does not have a classical (smooth) solution. A
notion of generalized solution-- viscosity solution has been developed to
cover such situations \textrm{\cite{cil,p3}}. This section will prove that
BSDEs under super/sub linear expectation provide a probabilistic
interpretation for the viscosity solution of HJB equations.

Consider the following backward stochastic differential equations under super
linear expectation coupled with a forward diffusion driven by a $d$%
-dimensional G-Brownian motion with components $(B_{t}^{j})$
\textbf{independent} of $(B_{t}^{i})$, $1\leq i<j\leq d$:%
\begin{align}
X_{s}^{t,x}  &  =x+\int_{t}^{s}b\left(  X_{r}^{t,x}\right)  dr+\int_{t}%
^{s}h_{j}\left(  X_{r}^{t,x}\right)  d\left\langle B^{j}\right\rangle
_{r}+\int_{t}^{s}\sigma_{j}\left(  X_{r}^{t,x}\right)  dB_{r}^{j}%
,\ t\in\lbrack0,T],\\
Y_{s}^{t,x}  &  =\mathbb{E}_{\ast}[\Phi(X_{T})+\int_{s}^{T}g\left(
X_{r}^{t,x},Y_{r}^{t,x}\right)  dr+\int_{s}^{T}f_{j}\left(  X_{r}^{t,x}%
,Y_{r}^{t,x}\right)  d\left\langle B^{j}\right\rangle _{r}|\mathcal{F}%
_{s}],s\in\lbrack t,T],
\end{align}
where $b,\sigma_{j},h_{j}:\mathbf{R}^{n}\mapsto\mathbf{R}^{n}$, $j=1,...d$,
$\Phi:\mathbf{R}^{n}\mapsto\mathbf{R}$, $g$, $f_{j}:\mathbf{R}^{n}%
\times\mathbf{R\mapsto R}$\textbf{. }Here and in the sequence we use the
Einstein convention, i.e., the above repeated indices $j$ implies summation.

We assume

{\textbf{(H4.1)}}. $\Phi(\cdot),b(\cdot),h_{j}(\cdot),\sigma_{j}%
(\cdot),g\left(  \cdot,\cdot\right)  ,f_{j}\left(  \cdot,\cdot\right)  $, are
given uniform Lipstchtiz functions. i.e., $\forall y^{1},y^{2}$:
$|\mathit{\phi}\left(  y^{1}\right)  -\mathit{\phi}\left(  y^{2}\right)  |\leq
K|y^{1}-y^{2}|$.

By Peng \textrm{\cite{p5,p9}} and Theorem 3.1 in this paper, there is a unique
pair $(X_{s}^{t,x},Y_{s}^{t,x})\in{\mathcal{M}_{G}^{\mathrm{2}}}%
(t,{\mathit{T}};\mathbf{R}^{n})\times{\mathcal{H}_{G}^{\mathrm{1}}%
}(t,{\mathit{T}};\mathbf{R})$ for (4.1) and (4.2). It is not difficult to
obtain the following estimates. Most of proofs can be founded in Peng
\textrm{\cite{p9}} Ch.V, Sec. 3.

\begin{lemma}
(i) $\mathbb{E}^{\ast}\left\vert X_{s}^{t,x}-X_{s}^{t,x^{\prime}}\right\vert
\leq C\left\vert x-x^{\prime}\right\vert $, $s\in\lbrack t,T]$.

(ii) $\mathbb{E}^{\ast}\left\vert X_{s}^{t,x}\right\vert ^{p}\leq
C(1+\left\vert x\right\vert ^{p})\delta^{\frac{p}{2}}$, $p\geq2$, $s\in\lbrack
t,T]$.

(iii) $\mathbb{E}^{\ast}\left\vert X_{t+\delta}^{t,x}-x\right\vert ^{2}\leq
C(1+\left\vert x\right\vert ^{2})\delta$, $\delta\in\lbrack0,T-t]$.

(iv) $\left\vert Y_{t}^{t,x}-Y_{t}^{t,x^{\prime}}\right\vert \leq C\left\vert
x-x^{\prime}\right\vert $.

(v) $\left\vert Y_{t}^{t,x}\right\vert \leq C(1+\left\vert x\right\vert )$.

(vi) $\left\vert Y_{t+\delta}^{t+\delta,x}-Y_{t}^{t,x}\right\vert \leq
C(1+\left\vert x\right\vert )(\delta^{\frac{1}{2}}+\delta)$.

(vii) \ $\mathbb{E}^{\ast}\left\vert Y_{t+\delta}^{t,x}-Y_{t}^{t,x}\right\vert
\leq C(1+\left\vert x\right\vert )\delta.$

where the constant C dose not depend on the $(t,x,\delta)$.
\end{lemma}

\bigskip Note that a nonlinear Feynman-Kac formula has been established in
\cite{p9}, where the control variable in associated PDE is a matrix. Here in
the framework of superlinear expectation, we assume some independence among
components of G-Brownian motion and obtain an obvious form of HJB equation
with vector-valued control variable. It is seen that in the next section
results here give some natural applications to stochastic control and
uncertainty volatility model.

We define
\begin{equation}
u\left(  t,x\right)  :=Y_{t}^{t,x}=\mathbb{E}_{\ast}\left[  \Phi(X_{T}%
^{t,x})+\int_{t}^{T}g\left(  X_{r}^{t,x},Y_{r}^{t,x}\right)  dr+\int_{t}%
^{T}f_{j}\left(  X_{r}^{t,x},Y_{r}^{t,x}\right)  d\left\langle B^{j}%
\right\rangle _{r}|\mathcal{F}_{t}\right]  .\nonumber
\end{equation}
Since $(X_{s}^{t,x},Y_{s}^{t,x})$ is dependant of $\mathcal{F}_{t}$, $u\left(
t,x\right)  $ is a deterministic function.

\begin{theorem}
$u\left(  t,x\right)  $ is a viscosity solution of the following HJB equation%
\begin{align}
\partial_{t}u+\underset{\alpha\in\Gamma}{\inf}\left\{  \mathcal{L}\left(
x,\alpha\right)  u+g\left(  x,u\right)  \right\}   &  =0\\
u|_{t=T}  &  =\Phi.
\end{align}
where the control variable $\alpha$ is selected dynamically within $\Gamma$.
$\mathcal{L}\left(  x,\alpha\right)  $ is a second order elliptic partial
differential operator parameterized by $\alpha$,
\[
\mathcal{L}\left(  x,\alpha\right)  =\frac{1}{2}\sum_{\mu,\nu=1}^{n}\left(
\sum_{j=1}^{d}\sigma_{\mu j}\sigma_{\nu j}\left(  x\right)  \alpha_{j}%
^{2}\right)  \cdot\partial_{x^{\mu}x^{\nu}}+\sum_{i=1}^{n}\left(  b_{i}\left(
x\right)  +\sum_{j=1}^{d}h_{ij}(x)\alpha_{j}^{2}\right)  \partial_{x^{i}}%
+\sum_{j=1}^{d}f_{j}\left(  x,u\right)  \alpha_{j}^{2}.
\]

\end{theorem}

\paragraph{Proof}

\bigskip(vi) and (iv) in Proposition 4.1 lead to the continuity of $u\left(
t,x\right)  $ in $\left(  t,x\right)  $. Now we prove $u\left(  t,x\right)  $
is a viscosity of PDE (4.3). By the definition of $u\left(  t,x\right)  $, for
$\delta\in(0,T-t)$, we have,%

\begin{align*}
u\left(  t+\delta,X_{t+\delta}^{t,x}\right)   &  =Y_{t+\delta}^{t+\delta
,X_{t+\delta}^{t,x}}=Y_{t+\delta}^{t,x}\\
&  =\mathbb{E}_{\ast}\left[  \Phi(X_{T}^{t,x})+\int_{t+\delta}^{T}g\left(
X_{r}^{t,x},Y_{r}^{t,x}\right)  dr+\int_{t+\delta}^{T}f_{j}\left(  X_{r}%
^{t,x},Y_{r}^{t,x}\right)  d\left\langle B^{j}\right\rangle _{r}%
|\mathcal{F}_{t+\delta}\right]  .
\end{align*}
Thus%
\[
u\left(  t,x\right)  =\mathbb{E}_{\ast}\left[  u\left(  t+\delta,X_{t+\delta
}^{t,x}\right)  +\int_{t}^{t+\delta}g\left(  X_{r}^{t,x},Y_{r}^{t,x}\right)
dr+\int_{t}^{t+\delta}f_{j}\left(  X_{r}^{t,x},Y_{r}^{t,x}\right)
d\left\langle B^{j}\right\rangle _{r}\right]  .
\]

\bigskip Now for fixed $\left(  t,x\right)  \in(0,T)\times\mathbf{R}^{n}$, Let
$\psi\in C_{b}^{1,3}([0,T]\times\mathbf{R}^{n})$ s.t $\psi\geq u$ and
$\psi\left(  t,x\right)  =u\left(  t,x\right)  $. By G-It\^{o} formula, it
follows that, for $\delta\in(0,T-t)$,%

\[
0\leq\mathbb{E}_{\ast}\left[  \psi\left(  t+\delta,X_{t+\delta}^{t,x}\right)
+\int_{t}^{t+\delta}g\left(  X_{r}^{t,x},Y_{r}^{t,x}\right)  dr+\int
_{t}^{t+\delta}f_{j}\left(  X_{r}^{t,x},Y_{r}^{t,x}\right)  d\left\langle
B^{j}\right\rangle _{r}-\psi\left(  t,x\right)  \right]  =\mathbb{E}_{\ast
}\left[  \mathbf{I}_{\delta}^{1}\right]
\]
Where we denote
\begin{align*}
\mathbf{I}_{\delta}^{1}  &  =\int_{t}^{t+\delta}\partial\psi_{t}\left(
r,X_{r}^{t,x}\right)  dr+\int_{t}^{t+\delta}\partial_{x^{i}}\psi\left(
r,X_{r}^{t,x}\right)  \left(  b_{i}\left(  X_{r}^{t,x}\right)  dr+h_{ij}%
\left(  X_{r}^{t,x}\right)  d\left\langle B^{j}\right\rangle _{r}\right) \\
&  +\frac{1}{2}\int_{t}^{t+\delta}\partial_{x^{\mu}x^{\nu}}\psi\left(
r,X_{r}^{t,x}\right)  \sigma_{\mu j}\sigma_{\nu j}\left(  X_{r}^{t,x}\right)
d\left\langle B\right\rangle _{r}^{j}\\
&  +\int_{t}^{t+\delta}g\left(  X_{r}^{t,x},Y_{r}^{t,x}\right)  dr+\int
_{t}^{t+\delta}f_{j}\left(  X_{r}^{t,x},Y_{r}^{t,x}\right)  d\left\langle
B^{j}\right\rangle _{r},
\end{align*}
and%
\begin{align*}
\mathbf{I}_{\delta}^{2}  &  =\delta\lbrack\partial\psi_{t}\left(  t,x\right)
+\partial_{x^{i}}\psi\left(  t,x\right)  b_{i}\left(  x\right)  +g\left(
x,u\left(  t,x\right)  \right)  ]\\
&  +[\frac{1}{2}\partial_{x^{\mu}x^{\nu}}\psi\left(  t,x\right)  \sigma_{\mu
j}\sigma_{\nu j}\left(  x\right)  +h_{ij}\left(  x\right)  +f_{j}\left(
x,u\left(  t,x\right)  \right)  ](\left\langle B^{j}\right\rangle _{t+\delta
}-\left\langle B^{j}\right\rangle _{t}).
\end{align*}
In the following we denote $C$ as a universal constant independent of $\delta
$. Then using the Lipschitz conditions and by Lemma 4.1 we have
\begin{align*}
\left\vert \mathbb{E}_{\ast}\left[  \mathbf{I}_{\delta}^{1}\right]
-\mathbb{E}_{\ast}\left[  \mathbf{I}_{\delta}^{2}\right]  \right\vert  &
\leq\mathbb{E}^{\ast}\left\vert \mathbf{I}_{\delta}^{1}-\mathbf{I}_{\delta
}^{2}\right\vert \\
&  \leq C\mathbb{E}^{\ast}[\int_{t}^{t+\delta}\left\vert Y_{r}^{t,x}%
-Y_{t}^{t,x}\right\vert dr+\int_{t}^{t+\delta}\left\vert X_{r}^{t,x}%
-x\right\vert dr\\
&  +\int_{t}^{t+\delta}(\left\vert X_{r}^{t,x}\right\vert +\left\vert
x\right\vert )\left\vert X_{r}^{t,x}-x\right\vert dr+\int_{t}^{t+\delta
}\left\vert X_{r}^{t,x}\right\vert ^{2}\left\vert X_{r}^{t,x}-x\right\vert
dr]\\
&  \leq C(1+\left\vert x\right\vert )\delta^{2}+C(1+\left\vert x\right\vert
+\left\vert x\right\vert ^{2}+\left\vert x\right\vert ^{3})\delta^{\frac{3}%
{2}}%
\end{align*}
We denote $R_{\delta}:=C(1+\left\vert x\right\vert )\delta+C(1+\left\vert
x\right\vert +\left\vert x\right\vert ^{2}+\left\vert x\right\vert ^{3}%
)\delta^{\frac{1}{2}}$. Then we obtain%

\begin{align*}
0  &  \leq\frac{1}{\delta}\mathbb{E}_{\ast}\left[  \mathbf{I}_{\delta}%
^{1}\right]  \leq\frac{1}{\delta}\mathbb{E}_{\ast}\left[  \mathbf{I}_{\delta
}^{3}\right]  +R_{\delta}\\
&  =\frac{1}{\delta}\mathbb{E}_{\ast}[\partial\psi_{t}\left(  t,x\right)
\delta+\partial_{x^{i}}\psi\left(  t,x\right)  \cdot b_{i}\left(  x\right)
\delta+g\left(  x,u\right)  \delta\\
&  +(\frac{1}{2}\partial_{x^{\mu}x^{\nu}}\psi\left(  r,x\right)  \sigma_{\mu
j}\sigma_{\nu j}\left(  x\right)  +\partial_{x^{i}}\psi\left(  t,x\right)
\cdot h_{ij}\left(  x\right)  +f_{j}\left(  x,u\right)  )\cdot(\left\langle
B^{j}\right\rangle _{t+\delta}-\left\langle B^{j}\right\rangle _{t}%
)]+R_{\delta}\\
&  =\partial_{t}u+\underset{\alpha\in\Gamma}{\inf}\left\{  \mathcal{L}\left(
x,\alpha\right)  u+g\left(  x,u\right)  \right\}  +R_{\delta}%
\end{align*}
Let $\delta\rightarrow0$, we deduce that $u\left(  t,x\right)  $ is a
viscosity subsolution of PDE (4.3). Similarly we can prove that $u\left(
t,x\right)  $ is also a viscosity supersolution of PDE (4.3). $\Box$

The HJB equation is the infinitesimal version of the dynamic programming
principle: it describes the local behavior of the super linear expectation.
The HJB equation is also called dynamic programming equation. Theorem 4.1
provides an immediate proof of the dynamic programming principle.

\begin{corollary}
Let $u$ be the solution of PDE (4.3). Then for every $\delta\in\lbrack0,T-t]$,
we have
\[
u\left(  t,x\right)  =\underset{P\in\mathcal{P}}{\inf}{\normalsize E}%
_{P}\left[  u\left(  t+\delta,X_{t+\delta}^{t,x}\right)  +\int_{t}^{t+\delta
}g\left(  X_{r}^{t,x},Y_{r}^{t,x}\right)  dr+\int_{t}^{t+\delta}f_{j}\left(
X_{r}^{t,x},Y_{r}^{t,x}\right)  d\left\langle B^{j}\right\rangle _{r}\right]
.
\]

\end{corollary}

\begin{corollary}
Let $(W_{t}),(B_{t}),(Z_{t})$ be three $G$-Brownian motion with uncertainty
set $\Gamma_{i},i=1,2,3$, with $(W_{t})$ independent of $(B_{t})$ and
$(Z_{t})$, $(B_{t})$ independent of $(Z_{t})$. Then the following system%
\begin{align}
X_{s}  &  =x+\int_{t}^{s}b\left(  X_{r}\right)  dr+\int_{t}^{s}h_{j}\left(
X_{r}\right)  d\left\langle W^{j}\right\rangle _{r}+\int_{t}^{s}\sigma
_{j}\left(  X_{r}\right)  dB_{r}^{j}\\
Y_{s}  &  =\mathbb{E}_{\ast}\left[  \Phi(X_{T})+\int_{s}^{T}g\left(
X_{r},Y_{r}\right)  dr+\int_{t}^{s}f_{j}\left(  X_{r},Y_{r}\right)
d\left\langle Z^{j}\right\rangle _{r}|\mathcal{F}_{s}\right]  ,
\end{align}
provide a probabilistic interpretation for the viscosity solution of following
equation%
\begin{align}
&  \partial_{t}u+\underset{\alpha\in\Gamma_{1}}{\inf}\left\{  \sum_{i=1}%
^{n}\sum_{j=1}^{d}h_{ij}(x)\alpha_{j}^{2}\partial_{x^{i}}u\right\}
+\underset{\alpha\in\Gamma_{2}}{\inf}\left\{  \frac{1}{2}\sum_{\mu,\nu=1}%
^{n}\left(  \sum_{j=1}^{d}\sigma_{\mu j}\sigma_{\nu j}\left(  x\right)
\alpha_{j}^{2}\right)  \partial_{x^{\mu}x^{\nu}}u\right\} \nonumber\\
&  \ \ \ \ +\underset{\alpha\in\Gamma_{3}}{\inf}\left\{  \sum_{j=1}^{d}%
f_{j}\left(  x,u\right)  \alpha_{j}^{2}\right\}  +\sum_{i=1}^{n}b_{i}\left(
x\right)  \partial_{x^{i}}u+g\left(  t,x,u\right)  =0\\
&  u|_{t=T}=\Phi.\
\end{align}

\end{corollary}

As a converse of Theorem 4.1, we have

\begin{theorem}
If PDE (4.3) has a classical solution $u\left(  t,x\right)  $ $\in C_{b}%
^{1,2}$, then $u\left(  s,X_{s}^{t,x}\right)  $ solves BSDE (4.2).
\end{theorem}

\paragraph{Proof}

We denote $N_{s}:=-\int_{t}^{s}\underset{\alpha\in\Gamma}{\inf}\left(
h_{ij}\partial_{x^{i}}u+\frac{1}{2}\partial_{x^{\mu}x^{\nu}}u\sigma_{\mu
j}\sigma_{\nu j}+f_{j}\right)  \alpha_{j}^{2}dr+\int_{t}^{s}(h_{ij}%
\partial_{x^{i}}u+\frac{1}{2}\partial_{x^{\mu}x^{\nu}}u\sigma_{\mu j}%
\sigma_{\nu j}+f_{j})d\left\langle B^{j}\right\rangle _{r}$. Here we use the
Einstein convention, i.e., the repeated indices $\mu,\nu,i,j$ imply summation.
Then by G-It\^{o} formula, we have%
\begin{align*}
u\left(  s,X_{s}^{t,x}\right)   &  =u\left(  t,x\right)  +\int_{t}^{s}\partial
u_{t}dr+\int_{t}^{s}\partial_{x^{i}}u\left(  b_{i}dr+h_{ij}d\left\langle
B^{j}\right\rangle _{r}+\sigma_{ij}dB_{r}^{j}\right)  +\frac{1}{2}\int_{t}%
^{s}\partial_{x^{\mu}x^{\nu}}u\sigma_{\mu j}\sigma_{\nu j}d\left\langle
B^{j}\right\rangle _{r}\\
&  =u\left(  t,x\right)  +\int_{t}^{s}\sigma_{ij}\partial_{x^{i}}udB_{r}%
^{j}+N_{s}-\int_{t}^{s}gdr-\int_{t}^{s}f_{j}d\left\langle B^{j}\right\rangle
_{r}%
\end{align*}
and%

\[
\Phi(X_{T})=u\left(  s,X_{T}^{t,x}\right)  =u\left(  t,x\right)  +\int_{t}%
^{T}\sigma_{ij}\partial_{x^{i}}udB_{r}^{j}+N_{T}-\int_{t}^{T}gdr-\int_{t}%
^{T}f_{j}d\left\langle B^{j}\right\rangle _{r}.
\]
Observe that $(N_{s})$ is a $\mathbb{E}_{\ast}$-martingale, thus
\begin{align*}
u\left(  s,X_{s}^{t,x}\right)   &  =\mathbb{E}_{\ast}\left[  \Phi(X_{T}%
)+\int_{t}^{T}g\left(  X_{r}^{t,x},u\left(  r,X_{r}^{t,x}\right)  \right)
dr+\int_{t}^{T}f_{j}\left(  X_{r}^{t,x},u\left(  r,X_{r}^{t,x}\right)
\right)  d\left\langle B^{j}\right\rangle _{r}|\mathcal{F}_{s}\right] \\
&  -\int_{t}^{s}g\left(  X_{r},Y_{r}\right)  dr-\int_{t}^{s}f_{j}\left(
X_{r},Y_{r}\right)  d\left\langle B^{j}\right\rangle _{r}\\
&  =\mathbb{E}_{\ast}\left[  \Phi(X_{T})+\int_{s}^{T}g\left(  X_{r}%
^{t,x},u\left(  r,X_{r}^{t,x}\right)  \right)  dr+\int_{s}^{T}f_{j}\left(
X_{r}^{t,x},u\left(  r,X_{r}^{t,x}\right)  \right)  d\left\langle
B^{j}\right\rangle _{r}|\mathcal{F}_{s}\right]  .
\end{align*}
Therefore $Y_{s}:=u\left(  s,X_{s}^{t,x}\right)  $ is the unique solution of
BSDE (4.2). $\Box$

\section{Applications}

With the Feynman-Kac formula, we now show some applications of
G-expectation(G-BSDE). In fact the condition of qusi-continuity is not
necessary for random variables and processes in sectin 3 and 4
\textrm{\cite{p9,lp}}. So we can apply all the results in section 3 and 4 to
practical problems not in qusi-continuous spaces.

\subsection{Connection with stochastic control problem}

Let $(\Omega,(\mathcal{F}_{t}),{\normalsize E})$ be a space of linear
expectation, $(W)_{t}$ is a $d$-dimensional standard ${\normalsize E}%
$-Brownian motion. $(\mathcal{F}_{t})_{t\geq0}$ is the usual augmented
Brownian motion filtration. Now consider the following stochastic control
problem with $n$-dimensional state process%
\begin{equation}
x_{s}=x+\int_{0}^{s}\left(  b\left(  x_{r}\right)  +h_{j}\left(  x_{r}\right)
\alpha_{j}^{2}(r)\right)  dr+\int_{0}^{s}\sigma_{j}\left(  x_{r}\right)
\alpha_{j}^{2}(r)dW_{r}^{j}%
\end{equation}
where $(\alpha_{s})_{s\geq0}$ $\in\mathcal{A}$, a set of $\mathcal{F}_{s}%
$-adapted processes taking values in a compact set $A\subseteq\mathbf{R}^{d}$,
called control process. $(x_{s})$ is called the trajectory corresponding to
$(\alpha_{s})$. Here we still use the Einstein convention. For any given
$t\in\lbrack0,T]$, we introduce the following one dimensional BSDE:%
\begin{equation}
y_{s}=E\left[  \Phi(x_{T-t})+\int_{s}^{T-t}\left(  g\left(  x_{r}%
,y_{r}\right)  +f_{j}\left(  x_{r},y_{r}\right)  \alpha_{j}^{2}(r)\right)
dr|\mathcal{F}_{s}\right]  ,\ s\in\lbrack0,T-t]
\end{equation}

\bigskip We assume

{\textbf{(H5.1)}}. $b,h$ and $\sigma$ are continuously differentiable in $x$,
their derivatives $b_{x}$, $h_{x}$, $\sigma_{x}$ being bounded.

{\textbf{(H5.2)}}. $f,g$ are continuously differentiable in $(x,y)$, their
derivatives $f_{x}$, $f_{y}$, $g_{x}$, $g_{y}$ being bounded.

\bigskip{\textbf{(H5.3)}}. $\Phi(x)$ is a uniform Lipschtiz function.

\bigskip Under conditions ({\textbf{H5.1)}}$\backsim${\textbf{(H5.3), }}the
system (5.1) and (5.2) is well defined and there is a unique pair
$(x_{s},y_{s})$ solving (5.1) and (5.2) (see \textrm{\cite{p3}}). Then we can
define the so called cost function as%

\[
J_{x,t}\left(  \alpha(\cdot)\right)  =y_{0}(={\normalsize E}y_{0})
\]
The value function of this optimal control problem is defined by%

\[
V\left(  x,T-t;\alpha(\cdot)\right)  =\underset{\alpha_{.}\in\mathcal{A}}%
{\inf}J\left(  x,t;\alpha(\cdot),\Phi(\cdot)\right)
\]
By Peng \textrm{\cite{p3}}, we have

\begin{proposition}
Let {\textbf{(H5.1)}}$\backsim${\textbf{(H5.3) }}hold. Then for fixed $\Phi$,
the value function $u(t,x):=V\left(  x,T-t;\Phi(\cdot)\right)  $,
$(t,x)\in\lbrack0,T]\times\mathbf{R}^{n}$ is a viscosity solution of HJB
equation (4.3) with control domain $A\subseteq\mathbf{R}^{d}$.
\end{proposition}

Then by the uniqueness of viscosity solution for PDE(4.3) (see Ishii
\textrm{\cite{ish}}), we have

\begin{theorem}
Let {\textbf{(H5.1)}}$\backsim${\textbf{(H5.3) }}hold and $A=\Gamma$ . Then
$V\left(  x,T-t;\Phi(\cdot)\right)  =u(t,x)=Y_{t}^{t,x}$.
\end{theorem}

A G-BSDE is in fact a recursive super (sub) strategy. Theorem 5.1 establishes
a direct connection between general recursive sub strategies and stochastic
control problems.

\subsection{The uncertainty volatility model}

\bigskip The uncertainty volatility model (UVM) for pricing and hedging
derivative securities in an environment where the volatility is not known
precisely, but is assumed instead to lie between two non-negative extreme
values $\underline{\sigma}$ and $\overline{\sigma}$. Let us denote
$\mathcal{P}$ the class of all probability measures ${\normalsize P}$ on the
set of path $\{S_{t},\ t\in\lbrack0,T]\}$ where $(S_{t})$ is the price of a
stock. For simplicity we restrict our discussion to derivative securities
based on a single liquidly traded stock which pays no dividends over the
contract's lifetime. Denote $\mathbb{E}_{\ast}\left[  \cdot\right]
=\underset{P\in\mathcal{P}}{\inf}{\normalsize E}_{P}\left[  \cdot\right]  $,
$\mathbb{E}^{\ast}\left[  \cdot\right]  =\underset{P\in\mathcal{P}}{\sup
}{\normalsize E}_{P}\left[  \cdot\right]  $. Let $(B_{t})_{t\geq0}$ be a
one-dimensional Brownian motion under $\mathbb{E}^{\ast}$ and $\mathbb{E}%
_{\ast}$. $(\mathcal{F}_{t})_{t\geq0}$ is the usual augmented Brownian motion
filtration. We denote parameters $\overline{\sigma}^{2}=\mathbb{E}^{\ast
}[\left\langle B\right\rangle _{t=1}]$, $\underline{\sigma}^{2}=\mathbb{E}%
_{\ast}[\left\langle B\right\rangle _{t=1}]$. If there is no arbitrage, the
forward stock price should satisfy the risk-neutral It\^{o} equation:%
\begin{align}
dS_{u}  &  =rS_{u}du+S_{u}dB_{u},u\in\lbrack t,T]\\
S_{t}  &  =x
\end{align}
where $r$ is the riskless interest rate.

Assume that at a given maturity date $T$, a derivative security is
characterized by $\Phi(S_{T})$, where $\Phi(\cdot)$ is a known function of the
price of the underlying stock. Then the offer price and bid price of this
derivative should satisfy%

\begin{equation}
\overline{Y_{s}}^{t,x}=\mathbb{E}^{\ast}\left[  \Phi(S_{T})+\int_{s}%
^{T}(-r\overline{Y_{u}}^{t,x})du|\mathcal{F}_{s}\right]  ,\ s\in\lbrack t,T]
\end{equation}

\begin{equation}
\underline{Y}_{s}^{t,x}=\mathbb{E}_{\ast}\left[  \Phi(S_{T})+\int_{s}%
^{T}(-r\underline{Y}_{u}^{t,x})du|\mathcal{F}_{s}\right]  ,\ s\in\lbrack t,T]
\end{equation}

However in Avellaneda, Levy and Par\'{a}s \textrm{\cite{alp} }or by Theorem
4.1, the offer price and bid price are characterized as $w^{\ast}(t,S_{t})$
and $w_{\ast}(t,S_{t})$ respectively:%

\begin{align}
\frac{\partial w^{\ast}}{\partial t}(t,x)+r(x\frac{\partial w^{\ast}}{\partial
x}(t,x)-w^{\ast}(t,x))+\frac{x^{2}}{2}\underset{\underline{\sigma}\leq
\sigma\leq\overline{\sigma}}{\sup}\left\{  \sigma^{{\normalsize 2}}%
\frac{\partial^{2}w^{\ast}}{\partial x^{2}}(t,x)\right\}   &  =0\\
w^{\ast}|_{t=T}  &  =\Phi.
\end{align}
and%

\begin{align}
\frac{\partial w_{\ast}}{\partial t}(t,x)+r(x\frac{\partial w_{\ast}}{\partial
x}(t,x)-w_{\ast}(t,x))+\frac{x^{2}}{2}\underset{\underline{\sigma}\leq
\sigma\leq\overline{\sigma}}{\inf}\left\{  \sigma^{{\normalsize 2}}%
\frac{\partial^{2}w_{\ast}}{\partial x^{2}}(t,x)\right\}   &  =0\\
w_{\ast}|_{t=T}  &  =\Phi.
\end{align}
which are referred as the Black-Scholes-Barrenblett (BSB) equations. They are
a generalization of the classical Black-Scholes PDE, and reduce to it in the
case of $\underline{\sigma}=\overline{\sigma}$.

By results in section 4, we have

\begin{proposition}
The two characterization for the offer price and bid price are equivalent,
i.e.:
\[
w^{\ast}(t,x)=\overline{Y_{t}}^{t,x}=\underset{P}{\sup}{\normalsize E}%
_{P}[e^{-r(T-t)}\Phi(X_{T})]
\]%
\[
w_{\ast}(t,x)=\underline{Y}_{t}^{t,x}=\underset{P}{\inf}{\normalsize E}%
_{P}[e^{-r(T-t)}\Phi(X_{T})]
\]

\end{proposition}

From the above we see that BSDEs under sub (super) linear expectation provide
a new characterization for the UVM model. Another problem is to solve either
BSDE(5.5) and (5.6) or BSB(5.7) and (5.9). When $\Phi$ is convex, the BSB
prices coincide the Black-Scholes prices at volatility $\overline{\sigma}$ and
$\underline{\sigma}$. Generally, we can not find explicit solutions for the
BSB equations. We now provide a representation for solutions of the
Black-Scholes-Barrenblett equations.

\begin{theorem}%
\[
w^{\ast}(t,x)=e^{-r(T-t)}\underset{\underline{\sigma}\leq\alpha_{.}%
\leq\overline{\sigma}}{\sup}{\normalsize E}\left[  \Phi
(xe^{r(T-t)+{\normalsize \int_{0}^{_{T-t}}}\alpha_{u}dW_{u}{\normalsize -}%
\frac{1}{2}{\normalsize \int_{0}^{s}}\left\vert {\normalsize \alpha_{u}%
}\right\vert ^{2}{\normalsize du}})\right]  ,
\]%
\[
w_{\ast}(t,x)=e^{-r(T-t)}\underset{\underline{\sigma}\leq\alpha_{.}%
\leq\overline{\sigma}}{\inf}{\normalsize E}\left[  \Phi
(xe^{r(T-t)+{\normalsize \int_{0}^{_{T-t}}}\alpha_{u}dW_{u}{\normalsize -}%
\frac{1}{2}{\normalsize \int_{0}^{s}}\left\vert {\normalsize \alpha_{u}%
}\right\vert ^{2}{\normalsize du}})\right]  .
\]

\end{theorem}

\paragraph{Proof}

Note that $x_{s}=x{\normalsize \exp\{rs+\int_{0}^{s}\alpha_{u}dW_{u}-}\frac
{1}{2}{\normalsize \int_{0}^{s}}\left\vert {\normalsize \alpha_{u}}\right\vert
^{2}{\normalsize du\}}$, $y_{s}=e^{-r(T-t)}E[\Phi(x_{T-t})|\mathcal{F}_{s}]$
solves%
\begin{align}
x_{s}  &  =x+\int_{0}^{s}rx_{u}du+\int_{0}^{s}\alpha_{u}x_{u}dW_{u},\\
y_{s}  &  ={\normalsize E}\left[  \Phi(x_{T-t})+\int_{s}^{T-t}(-ry_{u}%
)du|\mathcal{F}_{s}\right]  ,\ s\in\lbrack0,T-t],
\end{align}
where $(\alpha_{t})$ is the control process between $\underline{\sigma}$ and
$\overline{\sigma}$. $(W_{t})$ is a standard Brownian motion under linear
expectation ${\normalsize E}$. By Proposition 5.1, we have%

\begin{align*}
w_{\ast}(t,x)  &  =\underset{\underline{\sigma}\leq\alpha_{.}\leq
\overline{\sigma}}{\inf}{\normalsize E}y_{0}=e^{-r(T-t)}\underset
{\underline{\sigma}\leq\alpha_{.}\leq\overline{\sigma}}{\inf}{\normalsize E}%
\left[  \Phi(x_{T-t})\right] \\
&  =e^{-r(T-t)}\underset{\underline{\sigma}\leq\alpha_{.}\leq\overline{\sigma
}}{\inf}{\normalsize E}\left[  \Phi(xe^{r(T-t)+{\normalsize \int_{0}^{_{T-t}}%
}\alpha_{u}dW_{u}{\normalsize -}\frac{1}{2}{\normalsize \int_{0}^{s}%
}\left\vert {\normalsize \alpha_{u}}\right\vert ^{2}{\normalsize du}})\right]
\end{align*}
By Proposition 5.3, we have%

\begin{align*}
w^{\ast}(t,x)  &  =\overline{Y_{t}}^{t,x}=\underset{P}{\sup}{\normalsize E}%
_{p}[e^{-r(T-t)}\Phi(X_{T})]=-\underset{P}{\inf}{\normalsize E}_{p}%
[-e^{-r(T-t)}\Phi(X_{T})]=-\underline{Y}_{t}^{t,x}(-\Phi)\\
&  =e^{-r(T-t)}\underset{\underline{\sigma}\leq\alpha_{.}\leq\overline{\sigma
}}{\sup}{\normalsize E}\left[  \Phi(xe^{r(T-t)+{\normalsize \int_{0}^{_{T-t}}%
}\alpha_{u}dW_{u}{\normalsize -}\frac{1}{2}{\normalsize \int_{0}^{s}%
}\left\vert {\normalsize \alpha_{u}}\right\vert ^{2}{\normalsize du}})\right]
.
\end{align*}
$\Box$

\subsection{Representation for solutions of G-heat equations}

\bigskip We can use the following two PDEs%
\begin{equation}
\partial_{t}u(t,x)+G^{\ast}\left(  D_{x}^{2}u\right)  =0,\ u(T,x)=\Phi(x),
\label{5.17}%
\end{equation}%
\begin{equation}
\partial_{t}u(t,x)+G_{\ast}\left(  D_{x}^{2}u\right)  =0,\ u(T,x)=\Phi(x),
\label{5.18}%
\end{equation}
to introduce G-expectation (\textrm{\cite{p5}}), where $G^{\ast}\left(
a\right)  =\frac{1}{2}\underset{\underline{\sigma}\leq\sigma\leq
\overline{\sigma}}{\sup}\left\{  \sigma^{2}a\right\}  $, $G_{\ast}\left(
a\right)  =\frac{1}{2}\underset{\underline{\sigma}\leq\sigma\leq
\overline{\sigma}}{\inf}\left\{  \sigma^{2}a\right\}  $. So PDE(5.13) and
(5.14) are called G-heat equations. It is fundamentally important to solve
G-heat equations in the theory of G-expectation. Hu \textrm{\cite{hu}
}constructed explicit solutions of G$^{\ast}$-heat equation (5.13) with a
class of terminal condition $\Phi(x)=x^{n}$ for each integer $n\geq1$ and in
\textrm{\cite{dhp} }a representation for the solution of (5.13) is given on
the Wiener space. We now provide a representation for solutions of the above
G-heat equations (5.14) with terminal condition $\Phi\in C_{Lip}(\mathbf{R})$
under a given space of linear expectation. In fact results in this subsection
also holds for $\Phi\in C_{l,Lip}(\mathbf{R})$. The proof is a natural
sequence of Proposition 5.1 and Theorem 4.1. For simplicity we only consider
spatial variable $x\in\mathbf{R}$. Similarly to Theorem 5.2, we have

\begin{theorem}%
\[
u_{\ast}(t,x)=\underset{\underline{\sigma}\leq\alpha_{.}\leq\overline{\sigma}%
}{\inf}{\normalsize E}\left[  \Phi(x+\int_{0}^{T-t}\alpha_{r}dW_{r})\right]
\]
solves PDE(5.14).
\end{theorem}

\paragraph{Proof}

\bigskip Let $(W_{t})$ be a standard Brownian motion under a linear
expectation ${\normalsize E}$. Consider%
\begin{align*}
x_{s}  &  =x+\int_{0}^{s}\alpha_{r}dW_{r}\\
y_{s}  &  ={\normalsize E}\left[  \Phi(x_{T-t})|\mathcal{F}_{s}\right]
,\ s\in\lbrack0,T-t].
\end{align*}

By Proposition 5.1, we have%
\[
u_{\ast}(t,x)=\underset{\underline{\sigma}\leq\alpha_{.}\leq\overline{\sigma}%
}{\inf}{\normalsize E}\left[  \Phi(x_{T-t})\right]  =\underset{\underline
{\sigma}\leq\alpha_{.}\leq\overline{\sigma}}{\inf}{\normalsize E}\left[
\Phi(x+\int_{0}^{T-t}\alpha_{r}dW_{r})\right]  .
\]

\bigskip$\Box$

\bigskip Let $(B_{t})$ be a G-Brownian motion under $\mathbb{E}_{\ast}$, $X$
be a $\mathbb{E}_{\ast}$-normal distributed random variable.

\begin{corollary}
Assuming further $\Phi\in C^{2}(\mathbf{R})$, We have

(i)
\begin{align*}
\mathbb{E}_{\ast}[\Phi(x+(B_{T}-B_{t}))]  &  =\underset{\underline{\sigma}%
\leq\alpha_{.}\leq\overline{\sigma}}{\inf}{\normalsize E}\left[  \Phi
(x+\int_{0}^{T-t}\alpha_{r}dW_{r})\right] \\
&  =\Phi(x)+\frac{1}{2}\underset{\underline{\sigma}\leq\alpha_{.}\leq
\overline{\sigma}}{\inf}{\normalsize E}\left[  \int_{0}^{T-t}\Phi
^{\prime\prime}(x+\int_{0}^{r}\alpha_{u}dW_{u})\alpha_{r}^{2}dr\right]
\end{align*}

(ii)%
\begin{align*}
\mathbb{E}_{\ast}[\Phi(X))]  &  =\underset{\underline{\sigma}\leq\alpha
_{.}\leq\overline{\sigma}}{\inf}{\normalsize E}[\Phi(\int_{0}^{1}\alpha
_{r}dW_{r})]\\
&  =\Phi(0)+\frac{1}{2}\underset{\underline{\sigma}\leq\alpha_{.}\leq
\overline{\sigma}}{\inf}{\normalsize E}\left[  \int_{0}^{1}\Phi^{\prime\prime
}(\int_{0}^{r}\alpha_{u}dW_{u})\alpha_{r}^{2}dr\right]
\end{align*}

\end{corollary}

\paragraph{Proof}

\bigskip Consider%

\begin{align*}
X_{s}  &  =x+\int_{t}^{s}dB_{r}\\
Y_{s}  &  =\mathbb{E}_{\ast}\left[  \Phi(X_{T})|\mathcal{F}_{s}\right]
,\ s\in\lbrack t,T].
\end{align*}
By Theorem 5.1, we have $u(t,x)=Y_{t}^{t,x}=\mathbb{E}_{\ast}[\Phi
(X_{T})]=\mathbb{E}_{\ast}[\Phi((x+(B_{T}-B_{t}))]$ and $u(t,x)=V\left(
x,T-t;\Phi(\cdot)\right)  =\underset{\underline{\sigma}\leq\alpha_{.}%
\leq\overline{\sigma}}{\inf}{\normalsize E}[\Phi(x+\int_{0}^{T-t}\alpha
_{r}dW_{r})]$. So we obtain the desired results. $\Box$

\begin{remark}
When $\Phi$ is a convex function,%
\[
\mathbb{E}_{\ast}[\Phi(X))]=\frac{1}{\sqrt{2\pi}}\int_{-\infty}^{\infty}%
\Phi(\underline{\sigma}y)e^{-\frac{y^{2}}{2}}dy.
\]
and%
\[
\mathbb{E}^{\ast}[\Phi(X))]=\frac{1}{\sqrt{2\pi}}\int_{-\infty}^{\infty}%
\Phi(\overline{\sigma}y)e^{-\frac{y^{2}}{2}}dy,
\]
But if $\Phi$ is a concave function, $\underline{\sigma},\overline{\sigma}$
must interchange their positions. From the above we see that there is no fixed
density function for $G$-normal distributed random variable $X$ in the
traditional sense. Thus many problems under sub(super)-linear expectation can
not be solved via methods of density function.
\end{remark}

\section{Appendix: Discussion on the dominated convergence theorem}

The theorem of dominated convergence is fundamentally important in the theory
of classical probability. We initially want to embed a control process
$\alpha_{t}=\frac{d\left\langle B\right\rangle _{t}}{dt}$ into the coefficient
$g$ of BSDE (4.2). However when we derive the associated HJB equation, we find
that the dominated convergence theorem does not hold in general under
sublinear expectation induced by mutually singular probability measures.. We
now give a sufficient condition and a counterexample about it.

\begin{lemma}
Assume that a sequence $X_{n}$ converges to $X$ in the sense of `capacity',
i.e.: for any $\varepsilon>0$,
\[
C^{\ast}(\left\vert X_{n}-X\right\vert >\varepsilon)\rightarrow0.
\]
Then there is a subsequence $\{X_{n_{k}}\}$ such that $X_{n_{k}}\rightarrow
X,\ q.s.$.
\end{lemma}

\paragraph{Proof}

Since $X_{n}\rightarrow X$ in the sense of `capacity', we can find a
subsequence $\{X_{n_{k}}\}$ such that
\[
C^{\ast}(\left\vert X_{n_{k}}-X\right\vert >\frac{1}{k})\leq\frac{1}{k^{2}}.
\]
Then by the Borel-Cantelli Lemma (see Peng \textrm{\cite{p9}} Ch.VI, Lemma
1.5), we deduce that%
\[
C^{\ast}(\underset{k\rightarrow\infty}{\varlimsup}\{\left\vert X_{n_{k}%
}-X\right\vert >\frac{1}{k}\})=0.
\]
Obviously $X_{n_{k}}$ converges to $X$ on $(\underset{k\rightarrow\infty
}{\varlimsup}\{\left\vert X_{n_{k}}-X\right\vert >\frac{1}{k}\})^{c}$. $\Box$

Let $L{_{b}^{\mathrm{1}}}(\Omega)$ be the completion of all bounded
$\mathcal{B}(\Omega)-$measurable functions under norm $\Vert X\Vert
_{1}=\mathbb{E}^{\ast}\left\vert X\right\vert $

\begin{proposition}
Let $\{X_{n}\}_{n=1}^{\infty}$ be a measurable sequence on $(\Omega
,\mathcal{F})$ such that $\left\vert X_{n}\right\vert \leq Y$, q.s.,
n=1,2,\ldots\ and $Y\in L{_{b}^{\mathrm{1}}}(\Omega)$. If $X_{n}\rightarrow X$
in the sense of `capacity', then $\underset{n\rightarrow\infty}{\ \lim
}\mathbb{E}^{\ast}[X_{n}]=\mathbb{E}^{\ast}[X]$.
\end{proposition}

\paragraph{Proof}

By Lemma 6.1, there is a subsequence $\{X_{n_{k}}\}$ such that $X_{n_{k}%
}\rightarrow X,\ q.s.$. Therefore $\left\vert X\right\vert \leq Y$,
\textit{q.s.}\ and $X\in L{_{b}^{\mathrm{1}}}(\Omega)$. For any $\varepsilon
>0$, we denote $A_{n}:=\{\left\vert X_{n}-X\right\vert \geq\frac{\varepsilon
}{4}\}$. Then
\begin{equation}
\mathbb{E}^{\ast}[\left\vert X_{n}-X\right\vert \cdot\mathbf{1}%
_{\mathbf{\Omega/}A_{n}}]\leq\frac{\varepsilon}{4}C^{\ast}(\mathbf{\Omega
/}A_{n})<\frac{\varepsilon}{2}.
\end{equation}
By the complete continuity of $\mathbb{E}^{\ast}$ on $L{_{b}^{\mathrm{1}}%
}(\Omega)$ (see Peng \textrm{\cite{p9} }Ch.VI, Proposition 1.19), there exists
a $\delta>0$, for any $A\subset\Omega$, $C^{\ast}(A)<\delta$ such that
$\mathbb{E}^{\ast}[Y\mathbf{1}_{A}]<\frac{\varepsilon}{4}$. For this $\delta$,
there exists $n$, when $n\geq N$, $C^{\ast}(A_{n})<\delta$. Therefore%
\begin{equation}
\mathbb{E}^{\ast}[\left\vert X_{n}-X\right\vert \cdot\mathbf{1}_{A_{n}}%
]\leq2\mathbb{E}^{\ast}[Y\mathbf{1}_{A}]<\frac{\varepsilon}{2}.
\end{equation}
(6.1) and (6.2) lead to that%
\[
\left\vert \mathbb{E}^{\ast}X_{n}-\mathbb{E}^{\ast}X\right\vert \leq
\mathbb{E}^{\ast}[\left\vert X_{n}-X\right\vert \leq\mathbb{E}^{\ast
}[\left\vert X_{n}-X\right\vert \cdot\mathbf{1}_{\mathbf{\Omega/}A_{n}%
}]+\mathbb{E}^{\ast}[\left\vert X_{n}-X\right\vert \cdot\mathbf{1}_{A_{n}%
}]<\varepsilon,\ n\geq N.
\]
The proof is completed. $\Box$

`\textit{q.s}.' convergence does not implies `capacity' convergence. If we
replace `capacity' convergence by `\textit{q.s}.' convergence in Proposition
6.1, the dominated convergence theorem does not hold true.

\begin{counterexample}
Let $(B_{t})$ be a one dimensional G$_{[\underline{\sigma}^{2},\overline
{\sigma}^{2}]}$-Brownian motion with $\underline{\sigma}^{2}<\overline{\sigma
}^{2}$. For fixed $t$, we denote $X_{\delta}^{t}:=\frac{\left\langle
B\right\rangle _{t+\delta}-\left\langle B\right\rangle _{t}}{\delta}%
-\frac{\left\langle B\right\rangle _{t}-\left\langle B\right\rangle
_{t-\delta}}{\delta}$, $0<\delta<t$. Obviously $X_{\delta}^{t}\rightarrow
0,\ q.s.$ when $\delta\downarrow0$ and $\mathbb{E}^{\ast}[\underset
{\delta\downarrow0}{\ \lim}X_{\delta}^{t}]=0$. However since $\frac
{\left\langle B\right\rangle _{t+\delta}-\left\langle B\right\rangle _{t}%
}{\delta}$ is independent of $\frac{\left\langle B\right\rangle _{t}%
-\left\langle B^{{}}\right\rangle _{t-\delta}}{\delta}$, we have that
\begin{align*}
\underset{\delta\downarrow0}{\ \lim}\mathbb{E}^{\ast}[X_{\delta}^{t}] &
=\underset{\delta\downarrow0}{\ \lim}\mathbb{E}^{\ast}\left(  \frac
{\left\langle B\right\rangle _{t+\delta}-\left\langle B\right\rangle _{t}%
}{\delta}-\frac{\left\langle B\right\rangle _{t}-\left\langle B\right\rangle
_{t-\delta}}{\delta}\right)  \\
&  =\underset{\delta\downarrow0}{\ \lim}\mathbb{E}^{\ast}\left[
\mathbb{E}^{\ast}\left(  \frac{\left\langle B\right\rangle _{t+\delta
}-\left\langle B\right\rangle _{t}}{\delta}-y\right)  _{y=\frac{\left\langle
B\right\rangle _{t}-\left\langle B\right\rangle _{t-\delta}}{\delta}}\right]
\\
&  =\underset{\delta\downarrow0}{\ \lim}\mathbb{E}^{\ast}\left[
\overline{\sigma}^{2}-\frac{\left\langle B\right\rangle _{t}-\left\langle
B\right\rangle _{t-\delta}}{\delta}\right]  \\
&  =\overline{\sigma}^{2}-\mathbf{\underline{\sigma}}^{2}>0.
\end{align*}
From the above we see that $\underset{\delta\downarrow0}{\ \lim}%
\mathbb{E}^{\ast}[X_{\delta}^{t}]\neq\mathbb{E}^{\ast}[\underset
{\delta\downarrow0}{\ \lim}X_{\delta}^{t}]$. $\Box$
\end{counterexample}

\end{document}